\documentclass[10pt,a4paper]{article}

\usepackage{amsmath,amssymb,latexsym,pstricks,mathrsfs,comment,amsthm,graphicx,tikz,enumerate,
pgffor,cite,wrapfig,multicol,float}
\usepackage{bibspacing}
\newcommand{\nc}{\newcommand}
\nc{\rnc}{\renewcommand}

\usepackage{geometry}  \rnc{\ss}{\smallskip} \nc{\ms}{\medskip} \nc{\bs}{\bigskip} \nc{\nss}{\vspace{-3mm}}

\makeatletter

\begin{document}

\hyphenation{mon-oid mon-oids}

\nc{\Inj}{\operatorname{Inj}}
\nc{\Surj}{\operatorname{Surj}}
\nc{\U}{\mathcal U}
\nc{\trans}[2]{\left(\begin{smallmatrix} #1 \\ #2 \end{smallmatrix}\right)}
\nc{\bigtrans}[2]{\left(\begin{matrix} #1 \\ #2 \end{matrix}\right)}

\nc{\itemit}[1]{\item[\emph{(#1)}]}
\nc{\E}{\mathcal E}
\nc{\TX}{\T(X)}
\nc{\TXP}{\T(X,\P)}
\nc{\TXPi}{\T(X_i,\P_i)}
\nc{\EX}{\E(X)}
\nc{\EXP}{\E(X,\P)}
\nc{\EXPi}{\E(X_i,\P_i)}
\nc{\SX}{\S(X)}
\nc{\SXP}{\S(X,\P)}
\nc{\Sing}{\E}
\nc{\idrank}{\operatorname{idrank}}
\nc{\SingXP}{\Sing(X,\P)}
\nc{\De}{\Delta}
\nc{\sgp}{\operatorname{sgp}}
\nc{\mon}{\operatorname{mon}}
\nc{\Dn}{\mathcal D_n}
\nc{\Dm}{\mathcal D_m}

\nc{\lline}[1]{\draw(3*#1,0)--(3*#1+2,0);}
\nc{\uline}[1]{\draw(3*#1,5)--(3*#1+2,5);}
\nc{\thickline}[2]{\draw(3*#1,5)--(3*#2,0); \draw(3*#1+2,5)--(3*#2+2,0) ;}
\nc{\thicklabel}[3]{\draw(3*#1+1+3*#2*0.15-3*#1*0.15,4.25)node{{\tiny $#3$}};}

\nc{\slline}[3]{\draw(3*#1+#3,0+#2)--(3*#1+2+#3,0+#2);}
\nc{\suline}[3]{\draw(3*#1+#3,5+#2)--(3*#1+2+#3,5+#2);}
\nc{\sthickline}[4]{\draw(3*#1+#4,5+#3)--(3*#2+#4,0+#3); \draw(3*#1+2+#4,5+#3)--(3*#2+2+#4,0+#3) ;}
\nc{\sthicklabel}[5]{\draw(3*#1+1+3*#2*0.15-3*#1*0.15+#5,4.25+#4)node{{\tiny $#3$}};}

\nc{\stll}[5]{\sthickline{#1}{#2}{#4}{#5} \sthicklabel{#1}{#2}{#3}{#4}{#5}}
\nc{\tll}[3]{\stll{#1}{#2}{#3}00}

\nc{\ptranspic}[1]
{
\foreach \x in {0,2,4}
{\draw(0+\x,#1)--(1.75+\x,#1);}
}

\nc{\mfourpic}[9]{
\slline1{#9}0
\slline3{#9}0
\slline4{#9}0
\slline5{#9}0
\suline1{#9}0
\suline3{#9}0
\suline4{#9}0
\suline5{#9}0
\stll1{#1}{#5}{#9}{0}
\stll3{#2}{#6}{#9}{0}
\stll4{#3}{#7}{#9}{0}
\stll5{#4}{#8}{#9}{0}
\draw[dotted](6,0+#9)--(8,0+#9);
\draw[dotted](6,5+#9)--(8,5+#9);
}
\nc{\vdotted}[1]{
\draw[dotted](3*#1,10)--(3*#1,15);
\draw[dotted](3*#1+2,10)--(3*#1+2,15);
}
\nc{\vdottedx}[2]{
\draw[dotted](3*#1,#2*5)--(3*#1,#2*5+5);
\draw[dotted](3*#1+2,#2*5)--(3*#1+2,#2*5+5);
}

\nc{\mthreepic}[7]{
\slline1{#7}0
\slline3{#7}0
\slline4{#7}0
\suline1{#7}0
\suline3{#7}0
\suline4{#7}0
\stll1{#1}{#4}{#7}{0}
\stll3{#2}{#5}{#7}{0}
\stll4{#3}{#6}{#7}{0}
\draw[dotted](6,0+#7)--(8,0+#7);
\draw[dotted](6,5+#7)--(8,5+#7);
}

\nc{\Clab}[2]{
\sthicklabel{#1}{#1}{{}_{\phantom{#1}}C_{#1}}{1.25+5*#2}0
}
\nc{\sClab}[3]{
\sthicklabel{#1}{#1}{{}_{\phantom{#1}}C_{#1}}{1.25+5*#2}{#3}
}
\nc{\Clabl}[3]{
\sthicklabel{#1}{#1}{{}_{\phantom{#3}}C_{#3}}{1.25+5*#2}0
}
\nc{\sClabl}[4]{
\sthicklabel{#1}{#1}{{}_{\phantom{#4}}C_{#4}}{1.25+5*#2}{#3}
}
\nc{\Clabll}[3]{
\sthicklabel{#1}{#1}{C_{#3}}{1.25+5*#2}0
}
\nc{\sClabll}[4]{
\sthicklabel{#1}{#1}{C_{#3}}{1.25+5*#2}{#3}
}

\nc{\mtwopic}[6]{
\slline1{#6*5}{#5}
\slline2{#6*5}{#5}
\suline1{#6*5}{#5}
\suline2{#6*5}{#5}
\stll1{#1}{#3}{#6*5}{#5}
\stll2{#2}{#4}{#6*5}{#5}
}
\nc{\mtwopicl}[6]{
\slline1{#6*5}{#5}
\slline2{#6*5}{#5}
\suline1{#6*5}{#5}
\suline2{#6*5}{#5}
\stll1{#1}{#3}{#6*5}{#5}
\stll2{#2}{#4}{#6*5}{#5}
\sClabl1{#6}{#5}{i}
\sClabl2{#6}{#5}{j}
}

\nc{\keru}{\operatorname{ker}^\wedge} \nc{\kerl}{\operatorname{ker}_\vee}

\nc{\coker}{\operatorname{coker}}
\nc{\KER}{\ker}
\nc{\N}{\mathbb N}
\nc{\LaBn}{L_\al(\B_n)}
\nc{\RaBn}{R_\al(\B_n)}
\nc{\LaPBn}{L_\al(\PB_n)}
\nc{\RaPBn}{R_\al(\PB_n)}
\nc{\rhorBn}{\rho_r(\B_n)}
\nc{\DrBn}{D_r(\B_n)}
\nc{\DrPn}{D_r(\P_n)}
\nc{\DrPBn}{D_r(\PB_n)}
\nc{\DrKn}{D_r(\K_n)}
\nc{\alb}{\al_{\vee}}
\nc{\beb}{\be^{\wedge}}
\nc{\bnf}{\bn^\flat}
\nc{\Bal}{\operatorname{Bal}}
\nc{\Red}{\operatorname{Red}}
\nc{\Pnxi}{\P_n^\xi}
\nc{\Bnxi}{\B_n^\xi}
\nc{\PBnxi}{\PB_n^\xi}
\nc{\Knxi}{\K_n^\xi}
\nc{\C}{\mathscr C}
\nc{\exi}{e^\xi}
\nc{\Exi}{E^\xi}
\nc{\eximu}{e^\xi_\mu}
\nc{\Eximu}{E^\xi_\mu}
\nc{\REF}{ {\red [Ref?]} }
\nc{\GL}{\operatorname{GL}}
\rnc{\O}{\operatorname{O}}

\nc{\vtx}[2]{\fill (#1,#2)circle(.2);}
\nc{\lvtx}[2]{\fill (#1,0)circle(.2);}
\nc{\uvtx}[2]{\fill (#1,1.5)circle(.2);}

\nc{\Eq}{\mathfrak{Eq}}
\nc{\Gau}{\Ga^\wedge} \nc{\Gal}{\Ga_\vee}
\nc{\Lamu}{\Lam^\wedge} \nc{\Laml}{\Lam_\vee}
\nc{\bX}{{\bf X}}
\nc{\bY}{{\bf Y}}
\nc{\ds}{\displaystyle}

\nc{\uvert}[1]{\fill (#1,1.5)circle(.2);}
\nc{\uuvert}[1]{\fill (#1,3)circle(.2);}
\nc{\uuuvert}[1]{\fill (#1,4.5)circle(.2);}
\rnc{\lvert}[1]{\fill (#1,0)circle(.2);}
\nc{\overt}[1]{\fill (#1,0)circle(.1);}
\nc{\overtl}[3]{\node[vertex] (#3) at (#1,0) {  {\tiny $#2$} };}
\nc{\cv}[2]{\draw(#1,1.5) to [out=270,in=90] (#2,0);}
\nc{\cvs}[2]{\draw(#1,1.5) to [out=270+30,in=90+30] (#2,0);}
\nc{\ucv}[2]{\draw(#1,3) to [out=270,in=90] (#2,1.5);}
\nc{\uucv}[2]{\draw(#1,4.5) to [out=270,in=90] (#2,3);}
\nc{\textpartn}[1]{{\lower0.45 ex\hbox{\begin{tikzpicture}[xscale=.2,yscale=0.2] #1 \end{tikzpicture}}}}
\nc{\textpartnx}[2]{{\lower1.0 ex\hbox{\begin{tikzpicture}[xscale=.3,yscale=0.3] 
\foreach \x in {1,...,#1}
{ \uvert{\x} \lvert{\x} }
#2 \end{tikzpicture}}}}
\nc{\disppartnx}[2]{{\lower1.0 ex\hbox{\begin{tikzpicture}[scale=0.3] 
\foreach \x in {1,...,#1}
{ \uvert{\x} \lvert{\x} }
#2 \end{tikzpicture}}}}
\nc{\disppartnxd}[2]{{\lower2.1 ex\hbox{\begin{tikzpicture}[scale=0.3] 
\foreach \x in {1,...,#1}
{ \uuvert{\x} \uvert{\x} \lvert{\x} }
#2 \end{tikzpicture}}}}
\nc{\disppartnxdn}[2]{{\lower2.1 ex\hbox{\begin{tikzpicture}[scale=0.3] 
\foreach \x in {1,...,#1}
{ \uuvert{\x} \lvert{\x} }
#2 \end{tikzpicture}}}}
\nc{\disppartnxdd}[2]{{\lower3.6 ex\hbox{\begin{tikzpicture}[scale=0.3] 
\foreach \x in {1,...,#1}
{ \uuuvert{\x} \uuvert{\x} \uvert{\x} \lvert{\x} }
#2 \end{tikzpicture}}}}

\nc{\dispgax}[2]{{\lower0.0 ex\hbox{\begin{tikzpicture}[scale=0.3] 
#2
\foreach \x in {1,...,#1}
{\lvert{\x} }
 \end{tikzpicture}}}}
\nc{\textgax}[2]{{\lower0.4 ex\hbox{\begin{tikzpicture}[scale=0.3] 
#2
\foreach \x in {1,...,#1}
{\lvert{\x} }
 \end{tikzpicture}}}}
\nc{\textlinegraph}[2]{{\raise#1 ex\hbox{\begin{tikzpicture}[scale=0.8] 
#2
 \end{tikzpicture}}}}
\nc{\textlinegraphl}[2]{{\raise#1 ex\hbox{\begin{tikzpicture}[scale=0.8] 
\tikzstyle{vertex}=[circle,draw=black, fill=white, inner sep = 0.07cm]
#2
 \end{tikzpicture}}}}
\nc{\displinegraph}[1]{{\lower0.0 ex\hbox{\begin{tikzpicture}[scale=0.6] 
#1
 \end{tikzpicture}}}}
 
\nc{\disppartnthreeone}[1]{{\lower1.0 ex\hbox{\begin{tikzpicture}[scale=0.3] 
\foreach \x in {1,2,3,5,6}
{ \uvert{\x} }
\foreach \x in {1,2,4,5,6}
{ \lvert{\x} }
\draw[dotted] (3.5,1.5)--(4.5,1.5);
\draw[dotted] (2.5,0)--(3.5,0);
#1 \end{tikzpicture}}}}

\nc{\partn}[4]{\left( \begin{array}{c|c} 
#1 \ & \ #3 \ \ \\ \cline{2-2}
#2 \ & \ #4 \ \
\end{array} \!\!\! \right)}
\nc{\partnlong}[6]{\partn{#1}{#2}{#3,\ #4}{#5,\ #6}} 
\nc{\partnsh}[2]{\left( \begin{array}{c} 
#1 \\
#2 
\end{array} \right)}
\nc{\partncodefz}[3]{\partn{#1}{#2}{#3}{\emptyset}}
\nc{\partndefz}[3]{{\partn{#1}{#2}{\emptyset}{#3}}}
\nc{\partnlast}[2]{\left( \begin{array}{c|c}
#1 \ &  \ #2 \\
#1 \ &  \ #2
\end{array} \right)}

\nc{\uuarcx}[3]{\draw(#1,3)arc(180:270:#3) (#1+#3,3-#3)--(#2-#3,3-#3) (#2-#3,3-#3) arc(270:360:#3);}
\nc{\uuarc}[2]{\uuarcx{#1}{#2}{.4}}
\nc{\uuuarcx}[3]{\draw(#1,4.5)arc(180:270:#3) (#1+#3,4.5-#3)--(#2-#3,4.5-#3) (#2-#3,4.5-#3) arc(270:360:#3);}
\nc{\uuuarc}[2]{\uuuarcx{#1}{#2}{.4}}
\nc{\darcx}[3]{\draw(#1,0)arc(180:90:#3) (#1+#3,#3)--(#2-#3,#3) (#2-#3,#3) arc(90:0:#3);}
\nc{\darc}[2]{\darcx{#1}{#2}{.4}}
\nc{\udarcx}[3]{\draw(#1,1.5)arc(180:90:#3) (#1+#3,1.5+#3)--(#2-#3,1.5+#3) (#2-#3,1.5+#3) arc(90:0:#3);}
\nc{\udarc}[2]{\udarcx{#1}{#2}{.4}}
\nc{\uudarcx}[3]{\draw(#1,3)arc(180:90:#3) (#1+#3,3+#3)--(#2-#3,3+#3) (#2-#3,3+#3) arc(90:0:#3);}
\nc{\uudarc}[2]{\uudarcx{#1}{#2}{.4}}
\nc{\uarcx}[3]{\draw(#1,1.5)arc(180:270:#3) (#1+#3,1.5-#3)--(#2-#3,1.5-#3) (#2-#3,1.5-#3) arc(270:360:#3);}
\nc{\uarc}[2]{\uarcx{#1}{#2}{.4}}
\nc{\darcxhalf}[3]{\draw(#1,0)arc(180:90:#3) (#1+#3,#3)--(#2,#3) ;}
\nc{\darchalf}[2]{\darcxhalf{#1}{#2}{.4}}
\nc{\uarcxhalf}[3]{\draw(#1,1.5)arc(180:270:#3) (#1+#3,1.5-#3)--(#2,1.5-#3) ;}
\nc{\uarchalf}[2]{\uarcxhalf{#1}{#2}{.4}}
\nc{\uarcxhalfr}[3]{\draw (#1+#3,1.5-#3)--(#2-#3,1.5-#3) (#2-#3,1.5-#3) arc(270:360:#3);}
\nc{\uarchalfr}[2]{\uarcxhalfr{#1}{#2}{.4}}

\nc{\bdarcx}[3]{\draw[blue](#1,0)arc(180:90:#3) (#1+#3,#3)--(#2-#3,#3) (#2-#3,#3) arc(90:0:#3);}
\nc{\bdarc}[2]{\darcx{#1}{#2}{.4}}
\nc{\rduarcx}[3]{\draw[red](#1,0)arc(180:270:#3) (#1+#3,0-#3)--(#2-#3,0-#3) (#2-#3,0-#3) arc(270:360:#3);}
\nc{\rduarc}[2]{\uarcx{#1}{#2}{.4}}
\nc{\duarcx}[3]{\draw(#1,0)arc(180:270:#3) (#1+#3,0-#3)--(#2-#3,0-#3) (#2-#3,0-#3) arc(270:360:#3);}
\nc{\duarc}[2]{\uarcx{#1}{#2}{.4}}

\nc{\uv}[1]{\fill (#1,2)circle(.1);}
\nc{\lv}[1]{\fill (#1,0)circle(.1);}
\nc{\stline}[2]{\draw(#1,2)--(#2,0);}
\nc{\tlab}[2]{\draw(#1,2)node[above]{\tiny $#2$};}
\nc{\tudots}[1]{\draw(#1,2)node{$\cdots$};}
\nc{\tldots}[1]{\draw(#1,0)node{$\cdots$};}

\nc{\uvw}[1]{\fill[white] (#1,2)circle(.1);}
\nc{\huv}[1]{\fill (#1,1)circle(.1);}
\nc{\llv}[1]{\fill (#1,-2)circle(.1);}
\nc{\arcup}[2]{
\draw(#1,2)arc(180:270:.4) (#1+.4,1.6)--(#2-.4,1.6) (#2-.4,1.6) arc(270:360:.4);
}
\nc{\harcup}[2]{
\draw(#1,1)arc(180:270:.4) (#1+.4,.6)--(#2-.4,.6) (#2-.4,.6) arc(270:360:.4);
}
\nc{\arcdn}[2]{
\draw(#1,0)arc(180:90:.4) (#1+.4,.4)--(#2-.4,.4) (#2-.4,.4) arc(90:0:.4);
}
\nc{\cve}[2]{
\draw(#1,2) to [out=270,in=90] (#2,0);
}
\nc{\hcve}[2]{
\draw(#1,1) to [out=270,in=90] (#2,0);
}
\nc{\catarc}[3]{
\draw(#1,2)arc(180:270:#3) (#1+#3,2-#3)--(#2-#3,2-#3) (#2-#3,2-#3) arc(270:360:#3);
}

\nc{\arcr}[2]{
\draw[red](#1,0)arc(180:90:.4) (#1+.4,.4)--(#2-.4,.4) (#2-.4,.4) arc(90:0:.4);
}
\nc{\arcb}[2]{
\draw[blue](#1,2-2)arc(180:270:.4) (#1+.4,1.6-2)--(#2-.4,1.6-2) (#2-.4,1.6-2) arc(270:360:.4);
}
\nc{\loopr}[1]{
\draw[blue](#1,-2) edge [out=130,in=50,loop] ();
}
\nc{\loopb}[1]{
\draw[red](#1,-2) edge [out=180+130,in=180+50,loop] ();
}
\nc{\redto}[2]{\draw[red](#1,0)--(#2,0);}
\nc{\bluto}[2]{\draw[blue](#1,0)--(#2,0);}
\nc{\dotto}[2]{\draw[dotted](#1,0)--(#2,0);}
\nc{\lloopr}[2]{\draw[red](#1,0)arc(0:360:#2);}
\nc{\lloopb}[2]{\draw[blue](#1,0)arc(0:360:#2);}
\nc{\rloopr}[2]{\draw[red](#1,0)arc(-180:180:#2);}
\nc{\rloopb}[2]{\draw[blue](#1,0)arc(-180:180:#2);}
\nc{\uloopr}[2]{\draw[red](#1,0)arc(-270:270:#2);}
\nc{\uloopb}[2]{\draw[blue](#1,0)arc(-270:270:#2);}
\nc{\dloopr}[2]{\draw[red](#1,0)arc(-90:270:#2);}
\nc{\dloopb}[2]{\draw[blue](#1,0)arc(-90:270:#2);}
\nc{\llloopr}[2]{\draw[red](#1,0-2)arc(0:360:#2);}
\nc{\llloopb}[2]{\draw[blue](#1,0-2)arc(0:360:#2);}
\nc{\lrloopr}[2]{\draw[red](#1,0-2)arc(-180:180:#2);}
\nc{\lrloopb}[2]{\draw[blue](#1,0-2)arc(-180:180:#2);}
\nc{\ldloopr}[2]{\draw[red](#1,0-2)arc(-270:270:#2);}
\nc{\ldloopb}[2]{\draw[blue](#1,0-2)arc(-270:270:#2);}
\nc{\luloopr}[2]{\draw[red](#1,0-2)arc(-90:270:#2);}
\nc{\luloopb}[2]{\draw[blue](#1,0-2)arc(-90:270:#2);}

\nc{\larcb}[2]{
\draw[blue](#1,0-2)arc(180:90:.4) (#1+.4,.4-2)--(#2-.4,.4-2) (#2-.4,.4-2) arc(90:0:.4);
}
\nc{\larcr}[2]{
\draw[red](#1,2-2-2)arc(180:270:.4) (#1+.4,1.6-2-2)--(#2-.4,1.6-2-2) (#2-.4,1.6-2-2) arc(270:360:.4);
}

\rnc{\H}{\mathscr H}
\rnc{\L}{\mathscr L}
\nc{\R}{\mathscr R}
\nc{\D}{\mathcal D}
\nc{\J}{\mathscr D}

\nc{\ssim}{\mathrel{\raise0.25 ex\hbox{\oalign{$\approx$\crcr\noalign{\kern-0.84 ex}$\approx$}}}}
\nc{\POI}{\mathcal{POI}}
\nc{\wb}{\overline{w}}
\nc{\ub}{\overline{u}}
\nc{\vb}{\overline{v}}
\nc{\fb}{\overline{f}}
\nc{\gb}{\overline{g}}
\nc{\hb}{\overline{h}}
\nc{\pb}{\overline{p}}
\rnc{\sb}{\overline{s}}
\nc{\XR}{\pres{X}{R\,}}
\nc{\YQ}{\pres{Y}{Q}}
\nc{\ZP}{\pres{Z}{P\,}}
\nc{\XRone}{\pres{X_1}{R_1}}
\nc{\XRtwo}{\pres{X_2}{R_2}}
\nc{\XRthree}{\pres{X_1\cup X_2}{R_1\cup R_2\cup R_3}}
\nc{\er}{\eqref}
\nc{\larr}{\mathrel{\hspace{-0.35 ex}>\hspace{-1.1ex}-}\hspace{-0.35 ex}}
\nc{\rarr}{\mathrel{\hspace{-0.35 ex}-\hspace{-0.5ex}-\hspace{-2.3ex}>\hspace{-0.35 ex}}}
\nc{\lrarr}{\mathrel{\hspace{-0.35 ex}>\hspace{-1.1ex}-\hspace{-0.5ex}-\hspace{-2.3ex}>\hspace{-0.35 ex}}}
\nc{\nn}{\tag*{}}
\nc{\epfal}{\tag*{$\Box$}}
\nc{\tagd}[1]{\tag*{(#1)$'$}}
\nc{\ldb}{[\![}
\nc{\rdb}{]\!]}
\nc{\sm}{\setminus}
\nc{\I}{\mathcal I}
\nc{\InSn}{\I_n\setminus\S_n}
\nc{\dom}{\operatorname{dom}^\wedge} \nc{\codom}{\operatorname{dom}_\vee}
\nc{\ojin}{1\leq j<i\leq n}
\nc{\eh}{\widehat{e}}
\nc{\wh}{\widehat{w}}
\nc{\uh}{\widehat{u}}
\nc{\vh}{\widehat{v}}
\nc{\sh}{\widehat{s}}
\nc{\fh}{\widehat{f}}
\nc{\textres}[1]{\text{\emph{#1}}}
\nc{\aand}{\emph{\ and \ }}
\nc{\iif}{\emph{\ if \ }}
\nc{\textlarr}{\ \larr\ }
\nc{\textrarr}{\ \rarr\ }
\nc{\textlrarr}{\ \lrarr\ }

\nc{\comma}{,\ }

\nc{\COMMA}{,\qquad}
\nc{\TnSn}{\T_n\setminus\S_n} 
\nc{\TmSm}{\T_m\setminus\S_m} 
\nc{\TXSX}{\T_X\setminus\S_X} 
\rnc{\S}{\mathcal S}

\nc{\T}{\mathcal T} 
\nc{\A}{\mathscr A} 
\nc{\B}{\mathcal B} 
\rnc{\P}{\mathcal P} 
\nc{\K}{\mathcal K}
\nc{\PB}{\mathcal{PB}} 
\nc{\rank}{\operatorname{rank}}

\nc{\mtt}{\!\!\!\mt\!\!\!}

\nc{\sub}{\subseteq}
\nc{\la}{\langle}
\nc{\ra}{\rangle}
\nc{\mt}{\mapsto}
\nc{\im}{\mathrm{im}}
\nc{\id}{\mathrm{id}}
\nc{\bn}{[n]}
\nc{\ba}{[a]}
\nc{\bl}{[l]}
\nc{\bm}{[m]}
\nc{\bk}{[k]}
\nc{\br}{[r]}
\nc{\ve}{\varepsilon}
\nc{\al}{\alpha}
\nc{\be}{\beta}
\nc{\ga}{\gamma}
\nc{\Ga}{\Gamma}
\nc{\de}{\delta}
\nc{\ka}{\kappa}
\nc{\lam}{\lambda}
\nc{\Lam}{\Lambda}
\nc{\si}{\sigma}
\nc{\Si}{\Sigma}
\nc{\oijn}{1\leq i<j\leq n}
\nc{\oijm}{1\leq i<j\leq m}

\nc{\comm}{\rightleftharpoons}
\nc{\AND}{\qquad\text{and}\qquad}

\nc{\bit}{\vspace{-3 truemm}\begin{itemize}}
\nc{\smbit}{\vspace{-2 truemm}\begin{itemize}}
\nc{\bmc}{\vspace{-3 truemm}\begin{multicols}}
\nc{\emc}{\end{multicols}\vspace{-3 truemm}}
\nc{\eit}{\end{itemize}\vspace{-3 truemm}}
\nc{\ben}{\vspace{-3 truemm}\begin{enumerate}}
\nc{\een}{\end{enumerate}\vspace{-3 truemm}}
\nc{\eitres}{\end{itemize}}

\nc{\set}[2]{\{ {#1} : {#2} \}} 
\nc{\bigset}[2]{\big\{ {#1}: {#2} \big\}} 
\nc{\Bigset}[2]{\Big\{ \,{#1}\, \,\Big|\, \,{#2}\, \Big\}}

\nc{\pres}[2]{\la {#1} \,|\, {#2} \ra}
\nc{\bigpres}[2]{\big\la {#1} \,\big|\, {#2} \big\ra}
\nc{\Bigpres}[2]{\Big\la \,{#1}\, \,\Big|\, \,{#2}\, \Big\ra}
\nc{\Biggpres}[2]{\Bigg\la {#1} \,\Bigg|\, {#2} \Bigg\ra}

\nc{\pf}{\noindent{\bf Proof.}  }
\nc{\epf}{\hfill$\Box$\bigskip}
\nc{\epfres}{\hfill$\Box$}
\nc{\pfnb}{\pf}
\nc{\epfnb}{\bigskip}
\nc{\pfthm}[1]{\bigskip \noindent{\bf Proof of Theorem \ref{#1}}\,\,  } 
\nc{\pfprop}[1]{\bigskip \noindent{\bf Proof of Proposition \ref{#1}}\,\,  } 
\nc{\epfreseq}{\tag*{$\Box$}}

\makeatletter
\newcommand\footnoteref[1]{\protected@xdef\@thefnmark{\ref{#1}}\@footnotemark}
\makeatother

\numberwithin{equation}{section}

\newtheorem{thm}[equation]{Theorem}
\newtheorem{lemma}[equation]{Lemma}
\newtheorem{cor}[equation]{Corollary}
\newtheorem{prop}[equation]{Proposition}

\theoremstyle{definition}

\newtheorem{rem}[equation]{Remark}
\newtheorem{defn}[equation]{Definition}
\newtheorem{eg}[equation]{Example}
\newtheorem{ass}[equation]{Assumption}

\title{Idempotent rank in the endomorphism monoid of a non-uniform partition} 
\author{
Igor Dolinka%
\footnote{The first author gratefully acknowledges the support of Grant No.~174019 of the Ministry of Education, Science, and Technological Development of the Republic of Serbia, and Grant No.~1136/2014 of the Secretariat of Science and Technological Development of the Autonomous Province of Vojvodina.}
\\
{\footnotesize \emph{Department of Mathematics and Informatics}}\\
{\footnotesize \emph{University of Novi Sad, Trg Dositeja Obradovi\'ca 4, 21101 Novi Sad, Serbia}}\\
{\footnotesize {\tt dockie\,@\,dmi.uns.ac.rs}}\\~\\
James East%
\footnote{The second author gratefully acknowledges the support of the Glasgow Learning, Teaching, and Research Fund in partially funding his visit to the third author in July, 2014.}
\\
{\footnotesize \emph{Centre for Research in Mathematics; School of Computing, Engineering and Mathematics}}\\
{\footnotesize \emph{University of Western Sydney, Locked Bag 1797, Penrith NSW 2751, Australia}}\\
{\footnotesize {\tt J.East\,@\,uws.edu.au}}\\~\\
James D.~Mitchell\\
{\footnotesize \emph{Mathematical Institute, School of Mathematics and Statistics}}\\
{\footnotesize \emph{University of St Andrews, St Andrews, Fife KY16 9SS, United Kingdom}}\\
{\footnotesize {\tt jdm3\,@\,st-and.ac.uk}}
}

\maketitle

~\vspace{-10ex}
\begin{abstract}
We calculate the rank and idempotent rank of the semigroup $\mathcal E(X,\mathcal P)$ generated by the idempotents of the semigroup $\mathcal T(X,\mathcal P)$, which consists of all transformations of the finite set $X$ preserving a non-uniform partition $\mathcal P$.  We also classify and enumerate the idempotent generating sets of this minimal possible size.  This extends results of the first two authors in the uniform case.

%


\end{abstract}

\section{Introduction}\label{sect:intro}

%
Let $S$ be a monoid with identity $1$, and $E(S)=\set{s\in S}{s^2=s}$ the set of all idempotents of $S$.  For a subset $U\sub S$, we write $\la U\ra$ for the submonoid of $S$ generated by $U$, which consists of all products $u_1\cdots u_k$ with $u_1,\ldots,u_k\in U\cup\{1\}$. 
The \emph{rank} of $S$, denoted $\rank(S)$, is the minimal cardinality of a subset $U\sub S$ such that $S=\la U\ra$.  If $S$ is idempotent generated, then the \emph{idempotent rank} of $S$, denoted $\idrank(S)$, is the minimal cardinality of a subset $U\sub E(S)$ such that $S=\la U\ra$.  

The \emph{full transformation semigroup} on a set $X$, denoted $\T_X$, is the set of all transformations of $X$ (i.e., all functions $X\to X$), under the semigroup operation of composition.  
Let $\P=\set{C_i}{i\in I}$ be a partition of~$X$; that is, the sets $C_i$ are non-empty, pairwise disjoint, and their union is all of $X$.  The set
\[
\TXP = \set{f\in\T_X}{(\forall i\in I)(\exists j\in I) \ C_if\sub C_j},
\]
consisting of all transformations of $X$ preserving $\P$,
is a subsemigroup of $\T_X$.  
A calculation of $\rank(\TXP)$ for finite $X$ is given in \cite{AS2009} and \cite{ABMS2014} for the uniform and non-uniform cases, respectively.  (Recall that $\P$ is \emph{uniform} if $|C_i|=|C_j|$ for all $i,j\in I$.)  We write $\EXP=\la E(\TXP)\ra$ for the idempotent generated subsemigroup of $\TXP$.  In \cite{DE1}, the first two authors calculated $\rank(\EXP)$ and $\idrank(\EXP)$ in the case of $X$ being finite and $\P$ uniform; among other things, it was shown that the rank and idempotent rank are equal, and the idempotent generating sets of this minimal possible size were also classified and enumerated.  The purpose of the current work is to extend these results to the non-uniform case.  
Our main results include the classification and enumeration of the idempotents of $\TXP$ (Propositions~\ref{idempotent_prop} and~\ref{etxp_rec}); the calculation of the rank and idempotent rank of $\EXP$ (Theorem \ref{thm:rank} --- in particular, the rank and idempotent rank are equal unless $\P$ has exactly two blocks of size $1$ (and at least one other block)); and the classification and enumeration of all idempotent generating sets of the minimal possible size (Proposition \ref{prop:U} and Theorem~\ref{thm:igs_class}).  

\section{Preliminaries}\label{sect:TX}

In this section, we state a number of results we will need concerning $\T_X$ and $\TXP$ for uniform~$\P$.  For the remainder of the article, we fix a finite set $X$.  
The group of units of $\T_X$ is the \emph{symmetric group} $\S_X$, which consists of all permutations of $X$ (i.e., all bijections $X\to X$).  
Denote by $\E_X=\la E(\T_X)\ra$ the idempotent generated subsemigroup of $\T_X$.  We generally denote the identity element of any monoid by $1$; in particular, $1\in\T_X$ denotes the identity map on $X$, which we also sometimes write as $\id_X$.  If $x,y\in X$ and $x\not=y$, then we write $e_{xy}\in\T_X$ for the transformation defined by
\[
ze_{xy} = \begin{cases}
x &\text{if $z=y$}\\
z &\text{if $z\in X\sm\{y\}$.}
\end{cases}
\]
It is clear that $e_{xy}\in E(\T_X)$ for all $x,y$.  We write $\D_X=\set{e_{xy}}{x,y\in X,\;\!x\not=y}$.
The next result collects several facts from \cite{Howie1966,Howie1978,Gomes1987}.  We always interpret a binomial coefficient ${m\choose n}$ to be $0$ if $m<n$.  

\ms
\begin{thm}\label{thm:EX}
Let $X$ be a finite set with $|X|=n\geq0$.  Then
$$
\E_X= \la \D_X\ra = \{1\}\cup(\TXSX). 
$$
Further, $\rank(\E_X)=\idrank(\E_X)=\rho_n$, where $\rho_2=2$ and $\rho_n={n\choose2}$ if $n\not=2$. \epfres
\end{thm}

The minimal idempotent generating sets of $\E_X$ were characterised in \cite{Howie1978} in terms of strongly connected tournaments.  Such tournaments were enumerated in \cite{Wright1970}, and it was shown in \cite{DE1} that any idempotent generating set for $\E_X$ contains one of minimal size.

\ms
\begin{thm}\label{thm:sin}
Let $X$ be a finite set with $|X|=n\geq0$.  Then any idempotent generating set for $\E_X=\la E(\T_X)\ra$ contains an idempotent generating set of minimal possible size.  The number of minimal idempotent generating sets for $\E_X$ is equal to $\si_n$, where $\si_2=1$ and $\si_n=w_n$ for $n\not=2$, and where the numbers $w_n$ satisfy the recurrence
\[
w_0=1,\qquad w_n=F_n - \sum_{s=1}^{n-1}{n\choose s}w_{s}F_{n-s} \quad\text{for $n\geq1$},
\]
where $F_n=2^{{n\choose 2}}=2^{n(n-1)/2}$.  \epfres
\end{thm}

The following analogues of Theorems \ref{thm:EX} and \ref{thm:sin} were proved in \cite{DE1}.  

\ms
\begin{thm}\label{thm:rank_uniform}
Let $S=\EXP$, where $\P$ is a uniform partition of the finite set $X$ into $m\geq0$ blocks of size $n\geq1$.  Then $\rank(S)=\idrank(S)=\rho_{mn}$, where $\rho_{21}=2$ and $\rho_{mn}=m\rho_n + n!{m\choose2}$ if $(m,n)\not=(2,1)$.  The numbers $\rho_n$ are defined in Theorem \ref{thm:EX}. \epfres
\end{thm}

\ms
\begin{thm}\label{thm:enum_uniform}
Let $S=\EXP$, where $\P$ is a uniform partition of the finite set $X$ into $m\geq0$ blocks of size $n\geq1$.  Then any idempotent generating set of $S$ contains an idempotent generating set of minimal possible size.  The number of minimal idempotent generating sets for $S$ is equal to $\si_{mn}$, where
\[
\si_{mn} = \begin{cases}
1 &\text{if $m=0$}\\
\si_n &\text{if $m=1$}\\
\si_m &\text{if $n=1$}\\
\si_n^m\times \sum_{k=0}^{{m\choose 2}} w_{mk}(2^{n!}-2)^k &\text{if $m,n\geq2$.}
\end{cases}
\]
The numbers $\si_n$ are defined in Theorem \ref{thm:sin}, and the numbers $w_{nk}$ satisfy the recurrence 
\[
w_{00}=1,\qquad w_{nk} = F_{nk} - \sum_{s=1}^{n-1}{n\choose s} \sum_{l=0}^{k} w_{sl}F_{n-s,k-l} \quad\text{for $n\geq1$},
\]
where $F_{nk} = \displaystyle{{{n\choose2}\choose k}\cdot2^{{n\choose 2}-k}}$. \epfres
\end{thm}

\ms
\begin{rem}
It might seem odd to include the $m=0$ case (when $X=\emptyset$) in the previous two results, 
but these will be useful for stating and proving later results such as Theorems~\ref{thm:rank} and~\ref{thm:igs_class}. 
\end{rem}

\section{The semigroup $\EXP$}\label{sect:EXP}

For a non-negative integer $k$, we write $[k]=\{1,\ldots,k\}$, which we interpret to be empty if $k=0$.  We write $\T_k=\T_{[k]}$, and similarly for $\S_k$, $\E_k$, $\D_k$.  Denote by $\T(k,l)$ the set of all functions $[k]\to[l]$, noting that $\T(k,k)=\T_k$.  The \emph{image}, \emph{rank} and \emph{kernel} of a function $f:A\to B$ are defined by 
\[
\im(f) = \set{af}{a\in A} \COMMA
\rank(f) = |\im(f)| \COMMA 
\ker(f) = \set{(a,b)\in A\times A}{af=bf},
\]
respectively.  Obviously, 
\[
\im(fg)\sub\im(g) \COMMA \rank(fg)\leq\min(\rank(f),\rank(g)) \COMMA \ker(fg)\supseteq\ker(f)
\]
for all functions $f:A\to B$ and $g:B\to C$.  

Recall that $X$ is a fixed finite set.  We also fix a \emph{non-uniform} partition $\P=\{C_1,\ldots,C_m\}$ of $X$.  We will write $n_i=|C_i|$ for each $i$ and assume that $n_1\geq\cdots\geq n_m$.  We write $n=|X|=n_1+\cdots+n_m$.  
For convenience, we assume that $C_i=\{i\}\times[n_i]$ for each $i\in[m]$, so~$X=\set{(i,j)}{i\in[m],\;\!j\in[n_i]}$.

We now define some parameters associated to the partition $\P$ that will make statements of results cleaner (see especially Theorems \ref{thm:rank} and \ref{thm:igs_class}).  For $i\in[n]$, define the sets
\[
M_i=\set{q\in[m]}{n_q=i} \AND N_i=\set{j\in[i-1]}{M_j\not=\emptyset}=\set{n_q}{q\in[m],\;\!n_q<i},
\]
and put $\mu_i=|M_i|$ and $\nu_i=|N_i|$.  In particular, $\mu_i$ is the number of blocks of $\P$ of size $i$.  As an example, if $n=22$, $m=8$ and $(n_1,\ldots,n_8)=(5,5,3,2,2,2,2,1)$, then the values of $\mu_i,\nu_i$ are as follows:
\begin{center}
\begin{tabular}{c|cccccccccccccccccccccc}
$i$ & 1 & 2 & 3 & 4 & 5 & 6 & 7 & 8 & 9 & 10 &  11 & 12 & 13 & 14 & 15 & 16 & 17 & 18 & 19 & 20 &  21 & 22 \\
\hline
$\mu_i$ & 1 & 4 & 1 & 0 & 2 & 0 & 0 & 0 & 0 & 0 & 0 & 0 & 0 & 0 & 0 & 0 & 0 & 0 & 0 & 0 & 0 & 0 \\
$\nu_i$ & 0 & 1 & 2 & 3 & 3 & 4 & 4 & 4 & 4 & 4 & 4 & 4 & 4 & 4 & 4 & 4 & 4 & 4 & 4 & 4 & 4 & 4 \\
\end{tabular}
\end{center}


Let $f\in\TXP$.  There is a transformation $\fb\in\T_m$ such that, for all $i\in\bm$, $C_if\sub C_{i\fb}$.  Also, for each $i\in\bm$, there is a function $f_i\in\T(n_i,n_{i\fb})$ such that $(i,j)f=(i\fb,jf_i)$ for all $j\in[n_i]$.  The transformation $f\in\TXP$ is uniquely determined by $f_1,\ldots,f_m,\fb$, and we will write $f=[f_1,\ldots,f_m;\fb]$.  The product in $\TXP$ may easily be described in terms of this notation.  Indeed, if $f,g\in\TXP$, then $fg=[f_1g_{1\fb},\ldots,f_mg_{m\fb};\fb\gb]$.  
Note that $\overline{fg}=\fb\gb$ and $(fg)_i=f_ig_{i\fb}\in\T(n_i,n_{i\overline{fg}})$ for all $f,g\in\TXP$ and $i\in[m]$.  
When $\P$ is uniform, each $f_i$ belongs to $\T(n,n)=\T_n$ (where $n$ is the common size of each block of~$\P$), and $\TXP$ is a wreath product $\T_n\wr\T_m$, as noted in \cite{AS2009}.  

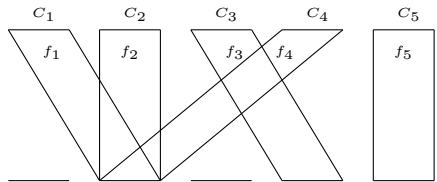
\begin{figure}[ht]
\begin{center}
\begin{tikzpicture}[xscale=.4,yscale=0.4]
\lline1
\lline2
\lline3
\lline4
\lline5
\uline1
\uline2
\uline3
\uline4
\uline5
\tll12{f_1}
\tll22{f_2}
\tll34{f_3}
\tll42{f_4}
\tll55{f_5}
\Clab10
\Clab20
\Clab30
\Clab40
\Clab50
\end{tikzpicture}
    \caption{Diagrammatic representation of an element of $\TXP$.}
    \label{fig:pic}
   \end{center}
 \end{figure}

\begin{figure}[ht]
\begin{center}
\begin{tikzpicture}[xscale=.4,yscale=0.4]
\lline1
\lline2
\lline3
\lline4
\lline5
\uline1
\uline2
\uline3
\uline4
\uline5
\tll11{g_1}
\tll23{g_2}
\tll31{g_3}
\tll44{g_4}
\tll54{g_5}
\slline1{2.5}{20}
\slline2{2.5}{20}
\slline3{2.5}{20}
\slline4{2.5}{20}
\slline5{2.5}{20}
\suline1{2.5}{20}
\suline2{2.5}{20}
\suline3{2.5}{20}
\suline4{2.5}{20}
\suline5{2.5}{20}
\stll13{f_1g_2}{2.5}{20}
\stll23{f_2g_2}{2.5}{20}
\stll34{f_3g_4}{2.5}{20}
\stll43{f_4g_2}{2.5}{20}
\stll54{f_5g_5}{2.5}{20}
\slline150
\slline250
\slline350
\slline450
\slline550
\suline150
\suline250
\suline350
\suline450
\suline550
\stll12{f_1}50
\stll22{f_2}50
\stll34{f_3}50
\stll42{f_4}50
\stll55{f_5}50
\draw(20,5)node{{\large $=$}};
\Clab11
\Clab21
\Clab31
\Clab41
\Clab51
\sClab1{.5}{20}
\sClab2{.5}{20}
\sClab3{.5}{20}
\sClab4{.5}{20}
\sClab5{.5}{20}
\end{tikzpicture}
    \caption{Diagrammatic calculation of a product in $\TXP$.}
    \label{fig:prod}
   \end{center}
 \end{figure}
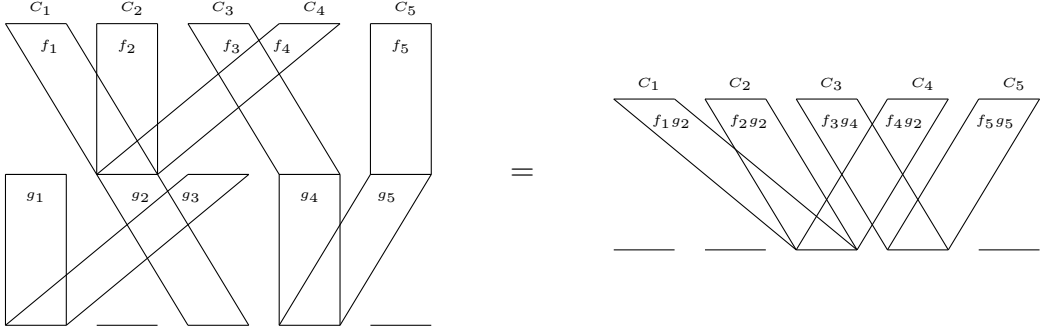

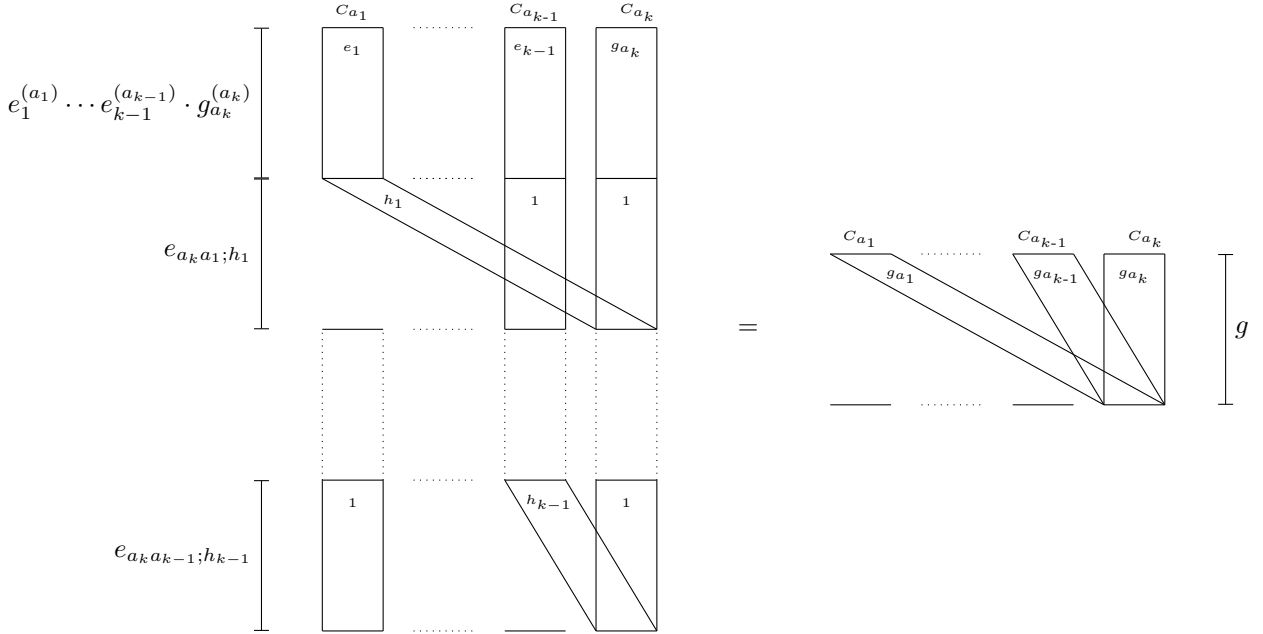
\begin{figure}[ht]
\begin{center}
\begin{tikzpicture}[scale=0.4]
\mthreepic134{e_1}{e_{k-1}}{g_{a_k}}{15}
\mthreepic434{h_1}11{10}
\mthreepic1441{h_{k-1}}1{0}
\vdottedx11
\vdottedx31
\vdottedx41
\draw(17,10)node{$=$};
\Clabll13{a_1}
\Clabll33{a_{k\text{-}1}}
\Clabl43{a_k}
\llbl0{e_{a_ka_{k-1};h_{k-1}}}
\llbl2{e_{a_ka_1;h_1}}
\llbl3{e_{1}^{(a_1)}\cdots e_{{k-1}}^{(a_{k-1})} \cdot g_{a_k}^{(a_k)}}
\end{tikzpicture}
\qquad
\begin{tikzpicture}[scale=0.4]
\mthreepic444{g_{a_1}}{g_{a_{k\text{-}1}}}{g_{a_k}}{15}
\draw[white](3,7.5)--(4,7.5); 
\Clabll1{3}{a_1}
\Clabll3{3}{a_{k\text{-}1}}
\Clabl4{3}{{a_k}}
\rlbl3{g}
\end{tikzpicture}
    \caption{Diagrammatic proof that $g=e_{1}^{(a_1)}\cdots e_{{k-1}}^{(a_{k-1})} \cdot g_{a_k}^{(a_k)} \cdot e_{a_ka_1;h_1}\cdots e_{a_ka_{k-1};h_{k-1}}$; see the proof of Proposition \ref{prop:G1G2} for more details.}
    \label{fig:h...he}
   \end{center}
 \end{figure}

There is a useful way to picture a transformation $f=[f_1,\ldots,f_m;\fb]\in\TXP$.  For example, with $m=5$, and $\fb=\left(\begin{smallmatrix} 1&2&3&4&5 \\ 2&2&4&2&5 \end{smallmatrix}\right)\in\T_5$, the transformation $f=[f_1,f_2,f_3,f_4,f_5;\fb]$ is pictured in Figure \ref{fig:pic}.  (Note that these diagrams are not supposed to imply that the sets $C_1,\ldots,C_m$ have the same size.)
This diagrammatic representation allows for easy visualisation of the multiplication.  For example, if $f$ is as above, and if $g=[g_1,g_2,g_3,g_4,g_5;\gb]$ where $\gb=\left(\begin{smallmatrix} 1&2&3&4&5 \\ 1&3&1&4&4 \end{smallmatrix}\right)$, then the product $fg=[f_1g_2,f_2g_2,f_3g_4,f_4g_2,f_5g_5;\fb\gb]$ may be calculated as in Figure \ref{fig:prod}.  Such diagrammatic methods may be used to verify various equations; an example is given in the proof of Proposition \ref{prop:G1G2} (see Figure \ref{fig:h...he}), but the rest are left to the reader.

The next result was proved in \cite[Proposition 3.1]{DE1} in the context of uniform partitions, but the argument works unmodified in the non-uniform case.

\ms
\begin{prop}\label{idempotent_prop}
A transformation $f\in\TXP$ is an idempotent if and only if
\bit
\itemit{i} $\fb\in E(\T_m)$,
\itemit{ii} $f_i\in E(\T_{n_i})$ for all $i\in\im(\fb)$, and
\itemit{iii} $\im(f_i) \sub \im(f_{i\fb})$ for all $i\in\bm\sm\im(\fb)$. \epfres
\eit
\end{prop}

A formula for $|E(\TXP)|$ was also given in \cite[Proposition 3.1]{DE1} in the uniform case.  That formula seems impossible to extend to the non-uniform case, but we may give a recurrence analogous to that of \cite[Proposition 3.2]{DE1}.  For a subset $A\sub[m]$, write $X_A=\bigcup_{a\in A}C_a$ and $\P_A=\set{C_a}{a\in A}$.  So~$\P_A$ is a partition of $X_A$ (which is empty if $A$ is empty).

\ms
\begin{prop}\label{etxp_rec}
Write $e(X,\P)=|E(\TXP)|$.  Then
\begin{align*}
e(X,\P) &= 1 &&\text{if $X$ is empty}\\
e(X,\P) &= \sum_{A} e(X_{A^c},\P_{A^c}) \sum_{a\in A} \sum_{l=1}^{n_{a}}{n_{a}\choose l}l^{n_A-l} &&\text{if $X$ is non-empty,}
\end{align*}
where the outer sum is over all $A\sub[m]$ with $1\in A$, and we write $A^c=[m]\sm A$ and $n_A=|X_A|=\sum_{a\in A}n_a$.
\end{prop}

\pf The statement for $X$ empty is clear, so suppose $X$ is non-empty.  An idempotent $f\in E(\TXP)$ is uniquely determined by:
\ben[(i)]
\item the set $A=\set{i\in\bm}{i\fb=1\fb}$,
\item the element $a=1\fb\in A$ (note that $1\fb\in A$ as $\fb$ is an idempotent),
\item the image $\im(f_{a})$, say of size $l\in[n_{a}]$ --- there are ${n_{a}\choose l}$ choices for these points, each of which are mapped identically by $f_{a}$,
\item the images under $f$ of the elements of $\big(\bigcup_{b\in A}C_b\big)\sm\big(\{a\}\times\im(f_{a})\big)$, which must all be in $\{a\}\times\im(f_{a})$ --- there are $l^{n_A-l}$ choices for these images, and then finally
\item the restriction of $f$ to $X_{A^c}=\bigcup_{i\in A^c}C_i$ --- this restriction belongs to $E(\T(X_{A^c},\P_{A^c}))$, which has size $e(X_{A^c},\P_{A^c})$.
\een
Multiplying these values and summing over relevant $A,a,l$ gives the desired result. \epf



We now move on to study the idempotent generated subsemigroup $\EXP=\la E(\TXP)\ra$ of $\TXP$.  
For simplicity, we will write $E=E(\TXP)$ and $S=\EXP=\la E\ra$.  

As in Section \ref{sect:TX}, for $k\geq2$ and $i,j\in[k]$ with $i\not=j$, we write $e_{ij}\in\T_k$ for the idempotent transformation defined by
\[
le_{ij} = \begin
{cases}
i &\text{if $l=j$}\\
l &\text{if $l\in[k]\sm\{j\}$.}
\end{cases}
\]
Note that $k$ (the size of the set on which $e_{ij}$ acts) depends on the context.
%
For non-negative integers $k,l$, we write $\Inj(k,l)$ (resp., $\Surj(k,l)$) for the set of all injective (resp., surjective) functions $[k]\to[l]$.  Note that $\Inj(k,k)=\Surj(k,k)=\S_k$, while if $k\not=l$, then (exactly) one of $\Inj(k,l)$ or $\Surj(k,l)$ is empty.

\begin{figure}[ht]
\begin{center}
\begin{tikzpicture}[xscale=.4,yscale=0.4]
\lline1
\lline2
\lline3
\lline4
\lline5
\uline1
\uline2
\uline3
\uline4
\uline5
\thickline11
\thickline24
\thickline33
\thickline44
\thickline55
\tll11{1}
\tll24{f}
\tll33{1}
\tll44{1}
\tll55{1}
\slline10{20}
\slline20{20}
\slline30{20}
\slline40{20}
\slline50{20}
\suline10{20}
\suline20{20}
\suline30{20}
\suline40{20}
\suline50{20}
\stll11{1}0{20}
\stll22{g}0{20}
\stll33{1}0{20}
\stll44{1}0{20}
\stll55{1}0{20}
\sClab10{20}
\sClab20{20}
\sClab30{20}
\sClab40{20}
\sClab50{20}
\Clab10
\Clab20
\Clab30
\Clab40
\Clab50
\end{tikzpicture}
    \caption{Diagrammatic representation of $e_{42;f}$ (left) and $g^{(2)}$ (right) from $\TXP$ with $m=5$.}
    \label{fig:e42f_f2}
   \end{center}
 \end{figure}
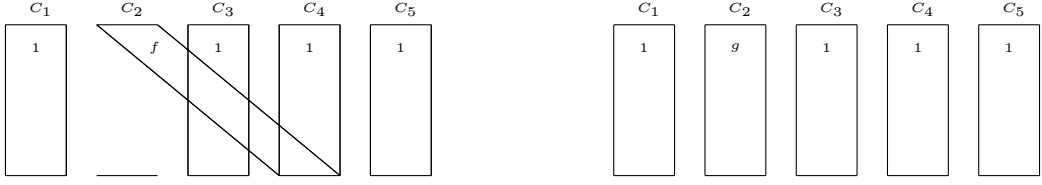

In what follows, certain special idempotents from $E$ will play a crucial role.  For $i,j\in\bm$ with $i\not=j$ and 
for any $f\in\Inj(n_j,n_i)$ or $\Surj(n_j,n_i)$, as appropriate, we write
\[
e_{ij;f} = [1,\ldots,1,f,1,\ldots,1;e_{ij}], 
\]
where $f$ is in the $j$th position. 
Note that here $e_{ij}=\overline{e_{ij;f}}$ refers to the idempotent $e_{ij}\in\T_m$.  The transformations $e_{ij;f}$ trivially satisfy conditions (i--iii) of Proposition \ref{idempotent_prop}, so $e_{ij;f}\in E$.  
If $k\in\bm$ and $g\in\T_{n_k}$, we will write $g^{(k)}=[1,\ldots,1,g,1,\ldots,1;1]$, where $g$ is in the $k$th position.  
For example, with $m=5$, the transformations $e_{42;f}$ and $g^{(2)}$ are pictured in Figure~\ref{fig:e42f_f2}.  
For any $k\in\bm$, and for any subset $U\sub \T_{n_k}$, we write $U^{(k)} = \set{g^{(k)}}{g\in U}$.  If $U$ is a subsemigroup of $\T_{n_k}$, then $U^{(k)}$ is a subsemigroup of $\TXP$ and is isomorphic to $U$.  Note that the $k$th coordinate of $e_{ij}^{(k)}=[1,\ldots,1,e_{ij},1,\ldots,1;1]$ is $e_{ij}\in\T_{n_k}$.

\ms
\begin{prop}\label{prop:G1G2}
The semigroup $S=\EXP$ is generated by $G_1\cup G_2$, where
\begin{align*}
G_1 = \set{e_{ij}^{(k)}}{k\in[m],\;\! i,j\in[n_k],\;\! i\not=j} \AND
G_2 = \set{e_{ij;f}}{i,j\in[m],\;\! i\not=j,\;\! f\in\Inj(n_j,n_i)\cup\Surj(n_j,n_i)}.
\end{align*}
\end{prop}

\pf Since the elements of $G_1\cup G_2$ are idempotents, it suffices to show that $E\sub\la G_1\cup G_2\ra$.  So let $f\in E$.  Write $A_1,\ldots,A_r$ for the $\ker(\fb)$-classes of $[m]$.  
Since $f$ is an idempotent, we have $f=f_1\cdots f_r$ where, for each $s\in[r]$, $f_s\in\TXP$ is defined by
\[
x f_s = \begin{cases}
xf &\text{if $x\in X_{A_s}=\bigcup_{a\in A_s}C_a$}\\
x &\text{if $x\in X\sm X_{A_s}$.}
\end{cases}
\]
So it suffices to show that $f_1,\ldots,f_r\in\la G_1\cup G_2\ra$.  Let $s\in[r]$, and write $A=A_s=\{a_1,\ldots,a_k\}$.  For simplicity, write $g=f_s=[g_1,\ldots,g_m;\gb]$.  Without loss of generality, suppose $A\gb=a_k$.  By Proposition~\ref{idempotent_prop} and Theorem~\ref{thm:EX}, we have $g_{a_k}\in E(\T_{n_{a_k}})\sub \la \D_{n_{a_k}}\ra$, and it quickly follows that $g_{a_k}^{(a_k)}\in\la G_1\ra$.  In particular, if $k=|A|=1$, then $g=g_{a_k}^{(a_k)}\in\la G_1\ra$.  So suppose $k\geq2$.  Now fix some $1\leq j<k$.  Let $e_{j}\in E(\T_{n_{a_j}})$ be such that $\ker(e_{j})=\ker(g_{a_j})$, and let $h_{j}\in\Inj(n_{a_j},n_{a_k})\cup\Surj(n_{a_j},n_{a_k})$ be any injective or surjective (as appropriate) map that extends the map $h_{j}':\im(e_{j})\to\im(g_{a_k})$ defined by $(xe_{j})h_{j}'=xg_{a_j}$ for $x\in[n_{a_j}]$.  In Figure~\ref{fig:h...he}, we show that
\[
g=e_{1}^{(a_1)}\cdots e_{{k-1}}^{(a_{k-1})} \cdot g_{a_k}^{(a_k)} \cdot e_{a_ka_1;h_1}\cdots e_{a_ka_{k-1};h_{k-1}}.
\]
(In the diagram, we only picture the action of the transformations on $X_A=C_{a_1}\cup\cdots\cup C_{a_k}$, and the pictured ordering of the blocks is not meant to imply that $a_1<\cdots<a_k$.)
Again, each $e_{j}^{(a_j)}$ belongs to $\la G_1\ra$, and clearly each $e_{a_ka_j;h_j}$ belongs to $G_2$.  This completes the proof. \epf

\nc{\llbl}[2]{
\draw[|-|] (1,#1*5)--(1,#1*5+5);
\draw(1,#1*5+2.5)node[left]{$#2$};
}
\nc{\rlbl}[2]{
\draw[|-|] (16,#1*5)--(16,#1*5+5);
\draw(16,#1*5+2.5)node[right]{$#2$};
}

Our next task is to calculate $\rank(S)$ and $\idrank(S)$.  In order to do this, we will show that the generating set $G_1\cup G_2$ from Proposition \ref{prop:G1G2} may be significantly reduced in size.  The next sequence of results (specifically, Lemmas~\ref{Elemma},~\ref{Glemma},~\ref{Glemma2} and~\ref{lem:U_r}) show what kind of transformations are essential in any (idempotent) generating set.  

Recall that $\E_{n_i}=\la E(\T_{n_i})\ra=\{1\}\cup(\T_{n_i}\sm\S_{n_i})$.  So $\E_{n_i}^{(i)}$, which consists of all maps $[1,\ldots,1,f,1,\ldots,1;1]\in \TXP$ with $f\in\E_{n_i}$ in the $i$th position, is a subsemigroup of $S$ isomorphic to $\E_{n_i}$.  The proof of \cite[Lemma 4.3]{DE1} is easily adapted to show the following.

\ms
\begin{lemma}\label{Elemma}
Let $i\in\bm$.  Then $S\sm\Sing_{n_i}^{(i)}$ is an ideal of $S$.  Consequently, any generating set for $S$ contains a generating set for $\E_{n_i}^{(i)}$. \epfres
\end{lemma}

Since the map $\TXP\to\T_m:f\mt\fb$ is a homomorphism, it follows from Proposition \ref{idempotent_prop} that $\fb\in\E_m=\{1\}\cup(\TmSm)$ for all $f\in S$.  We will frequently make use of this fact.  The next simple result describes the preimage of $1\in\T_m$ under the above map.

\ms
\begin{lemma}\label{lem:fb=1}
Let $f\in S$.  If $\fb=1$, then $f_i\in\E_{n_i}$ for all $i\in[m]$.
\end{lemma}

\pf Let $f=h_1\cdots h_k$, where $h_1,\ldots,h_k\in E$, and write $h_j=[h_{j1},\ldots,h_{jm};\hb_j]$ for each $j$.  Since $1=\fb=\hb_1\cdots\hb_k$, we see that $\hb_j=1$ for all $j$.  It follows that $h_{ji}\in E(\T_{n_i})$ for each $i,j$.  So $f_i=h_{1i}\cdots h_{ki}\in\la E(\T_{n_i})\ra=\E_{n_i}$ for each $i$.~\epf

%
%

For $\oijm$, we write $\ve_{ij}=\ve_{ji}$ for the equivalence relation on $\bm$ with unique non-trivial equivalence class $\{i,j\}$.  Note that $\ker(e_{ij})=\ker(e_{ji})=\ve_{ij}$.  We also write $\De = \set{(i,i)}{i\in\bm}$ for the trivial equivalence on $\bm$ (i.e., the equality relation on $\bm$).  


\ms
\begin{lemma}\label{Glemma}
Let $\oijm$ and $f\in\Inj(n_j,n_i)$, and suppose $e_{ij;f}=gh$ where $g,h\in S$ and $g\not=1$.  Then
\smbit\begin{multicols}{3}
\itemit{i} $\ker(\gb)=\ve_{ij}$,
\itemit{ii} $g_1,\ldots,g_m$ are injective, 
\itemit{iii} $g_i\in\S_{n_i}$ and $g_jg_i^{-1}=f$.
\end{multicols}\eit
Consequently, any generating set for $S$ contains such an element $g$ for each such $i,j,f$.  Further, if $g$ is an idempotent, then
\bit
\itemit{iv} if $n_i=n_j$, then either $g=e_{ij;f}$ or $g=e_{ji;f^{-1}}$,
\itemit{v} if $n_i>n_j$, then $g=e_{ij;f}$.
\eit
\end{lemma}

\pf 
Now,
$
[1,\ldots,1,f,1,\ldots,1;e_{ij}] = e_{ij;f} = gh = [g_1h_{1\gb},\ldots,g_mh_{m\gb};\gb\hb].
$
Since each $g_rh_{r\gb}$ is injective (equal to either $1$ or $f$), it follows that each $g_r$ is injective, establishing (ii).  If $\gb=1$, then $g_r\in\E_{n_r}$ for each~$r$ by Lemma~\ref{lem:fb=1}; but the only injective element of $\E_{n_r}$ is the identity element~$1$, so $g_r=1$ for all $r$, giving $g=[1,\ldots,1;1]=1$, a contradiction.  So $\gb\not=1$, whence $\gb\in\TmSm$.  But then
$
\De \not= \ker(\gb) \sub \ker(\gb \hb) = \ker(e_{ij}) = \ve_{ij},
$
so that $\ker(\gb)=\ve_{ij}$, giving (i).  Put $k=i\gb=j\gb$.

Next we claim that $n_k=n_i$.  Indeed, suppose this was not the case.  Since $g_i\in\Inj(n_i,n_k)$, we have $n_i\leq n_k$, so we must in fact have $n_i<n_k$.  Let $L=\set{r\in[m]}{n_r>n_i}$, noting that $k\in L$.  Since $n_1\geq\cdots\geq n_m$, it follows that $L=[s]$ for some $s\geq1$.  Now, since $\gb$ maps $[m]\sm\{i,j\}$ injectively into $[m]\sm\{k\}$, and since $s<i<j$, it follows that there exists $r\in L=[s]$ such that $r\gb>s$.  But then $g_r\in\T(n_r,n_{r\gb})$ with $n_{r\gb}\leq n_i<n_r$, contradicting the fact that $g_r$ is injective.  This completes the proof of the claim.

In particular, it follows that $g_i\in\Inj(n_i,n_k)=\Inj(n_i,n_i)=\S_{n_i}$.  Also, $h_k=h_{i\gb}\in\T(n_k,n_i)$ since $i=ie_{ij}=i\gb\hb=k\hb$.  We also have $1=g_ih_k$, so that $h_k=g_i^{-1}$, from which it follows that
$f=g_jh_{j\gb} = g_jh_k = g_jg_i^{-1}$, completing the proof of (iii).

Next, suppose $G$ is an arbitrary generating set for $S$.  By considering an expression $e_{ij;f}=h_1\cdots h_k$, where $h_1,\ldots,h_k\in G\sm\{1\}$, we see that $h_1\in G$ satisfies conditions (i--iii).  

Finally, suppose $g$ is an idempotent.  Since $\gb\in E(\T_m)$ by Proposition \ref{idempotent_prop}, and since $\ker(\gb)=\ve_{ij}$, it follows that $\gb=e_{ij}$ or $\gb=e_{ji}$.  Suppose first that $\gb=e_{ij}$.  Proposition \ref{idempotent_prop} also gives $g_r\in E(\T_{n_r})$ for each $r\in\im(\gb)=[m]\sm\{j\}$; but each $g_r$ is injective, so it follows that $g_r=1$ if $r\not=j$.  We also have $f=g_jg_i^{-1}=g_j$, giving $g=e_{ij;f}$.  Next suppose $\gb=e_{ji}$.  In particular, $n_i=n_j$ as $g_i\in\Inj(n_i,n_j)$ and $n_i\geq n_j$.  Again, we have $g_r=1$ for all $r\in[m]\sm\{i\}$, and this time we have $f=g_jg_i^{-1}=g_i^{-1}$, which gives $g_i=f^{-1}$, and $g=e_{ji;f^{-1}}$.  This completes the proof. \epf

Let $i,j\in[m]$ with $n_i>n_j$.  As a consequence of the previous result, we see that any idempotent generating set $G$ for $S$ contains all of
$
\set{e_{ik;f}}{k\in[m],\;\!n_k=n_j,\;\!f\in\Inj(n_j,n_i)}.
$
It might be tempting to guess that $G$ must also contain all of
$
\set{e_{ki;f}}{k\in[m],\;\!n_k=n_j,\;\!f\in\Surj(n_i,n_j)}.
$
But this is far from the case.  In fact, $G$ need only contain a \emph{single} element of the latter set, as we show in Lemma \ref{Glemma2} and Proposition \ref{prop:U}, the proof of which requires the next technical result (which will also be useful elsewhere).


\ms
\begin{lemma}\label{lem:f|L}
Let $j\in[n]$ and put $L=\set{q\in[m]}{n_q>j}$.  Suppose $f\in S$ is such that $\fb|_L$ is injective and $|C_qf|>j$ for all $q\in L$.  Then $\fb|_L=\id_L$ and $f_q\in\E_{n_q}$ for all $q\in L$.
\end{lemma}

\pf If $L=\emptyset$, there is nothing to show, so suppose $L\not=\emptyset$.  By Proposition \ref{prop:G1G2}, we may write $f=g_1\cdots g_k$ where  $g_1,\ldots,g_k\in E$ and $\rank(\gb_l)\geq m-1$ for each $l\in[k]$.  
(We could be more specific by insisting that $g_1,\ldots,g_k\in G_1\cup G_2$, but it will be convenient later to argue more generally, as we do here.)
We claim that $\gb_l|_L=\id_L$ for each $l$.  Indeed, suppose this is not the case, and let $l\in[k]$ be minimal so that $\gb_l|_L\not=\id_L$.  
Since $\gb_l\not=1$, it follows that $\rank(\gb_l)=m-1$, so $\gb_l=e_{ab}$ for some $a,b\in[m]$ with $a\not=b$.  
Note that at least one of $a,b$ belongs to $L$ since $\gb_l|_L\not=\id_L$.
We could not have $a,b\in L$ or else then $(a\gb_1\cdots\gb_{l-1})\gb_l=a\gb_l=b\gb_l=(b\gb_1\cdots\gb_{l-1})\gb_l$, giving $a\fb=b\fb$, contradicting the assumption that $\fb|_L$ is injective.  We also could not have $a\in L$ or else then $b\in[m]\sm L$ so that $\gb_l|_L=\id_L$, another contradiction.  So $a\in[m]\sm L$ and $b\in L$.  But then 
$|C_bf|=|C_bg_1\cdots g_k|\leq|(C_bg_1\cdots g_{l-1})g_l|\leq|C_bg_l|\leq|C_a|\leq j$, contradicting the assumption that $|C_qf|> j$ for all $q\in L$.  This establishes the claim that $\gb_l|_L=\id_L$ for each $l$.  It obviously follows that $\fb|_L=\id_L$.  For $l\in[k]$, write $g_l=[g_{l1},\ldots,g_{lm};\gb_l]$.  
Since $L\sub\im(\gb_l)$ and $g_l\in E$ for each $l\in[k]$, it follows from Proposition \ref{idempotent_prop} that $g_{lq}\in E(\T_{n_q})$ for each $l\in[k]$ and $q\in L$.
Since $\gb_l|_L=\id_L$ for all $l\in[k]$, it follows that $f_q=g_{1q}\cdots g_{kq}\in\E_{n_q}$ for each $q\in L$.  \epf


\begin{lemma}\label{Glemma2}
Let $i,j\in[m]$ be such that $n_i>n_j$, and let $f\in\Surj(n_i,n_j)$.  Suppose $e_{ji;f}=g_1\cdots g_r$ where $g_1,\ldots,g_r\in S\sm\{1\}$.  Put $L=\set{q\in[m]}{n_q>n_j}$.  Let $l\in[r]$ be minimal so that $\gb_l|_L\not=\id_L$.  Let $h=g_l$ and write $h=[h_1,\ldots,h_m;\hb]$.  Then 
\ms\begin{itemize}\bmc3
\itemit{i} $\rank(\hb)=m-1$, 
\itemit{ii} $|C_ih|=n_j$, and
\itemit{iii} $h_q$ is injective for all $q\in L\sm \{i\}$.
\emc\eit
%
Consequently, any generating set for $S$ contains such an element $h$ for each such $i,j$.  Further, if $h$ is an idempotent, then 
$h=e_{ki;f'}$ for some $k\in[m]$ with $n_k=n_j$ and some $f'\in\Surj(n_i,n_j)$.
\end{lemma}

\pf First note that $\hb\not=1$ since $\hb|_L\not=\id_L$, so $\hb\in\TmSm$ and $\rank(\hb)\leq m-1$.  But also $\rank(\hb)\geq\rank(e_{ji})=m-1$, so (i) holds.

Next, suppose (ii) does not hold.  First note that $n_j=|C_j|=|C_ie_{ji;f}|=|C_ig_1\cdots g_r|\leq|(C_ig_1\cdots g_{l-1})h|\leq|C_ih|$, since $\gb_1\cdots\gb_{l-1}$ acts as the identity on $L$ and $i\in L$.  So, since we are assuming that $|C_ih|\not=n_j$, we must have $|C_ih|>n_j$.  A similar calculation shows that $n_q\leq|C_qh|$ for any $q\in L\sm\{i\}$.  In particular, together with the assumption that $|C_ih|>n_j$, this gives $|C_qh|>n_j$ for all $q\in L$.  Since $(\gb_1\cdots\gb_{l-1}\hb\gb_{l+1}\cdots\gb_r)|_L=e_{ji}|_L$ is injective, and since $(\gb_1\cdots\gb_{l-1})|_L=\id_L$, it follows that $\hb|_L$ is injective.  But then Lemma \ref{lem:f|L} says that $\hb|_L=\id_L$, a contradiction.  This completes the proof of (ii).

Next, let $q\in L\sm\{i\}$ be arbitrary.  We have already seen that $|C_qh|\geq n_q$.  But we also trivially have $|C_qh|\leq |C_q|=n_q$, whence $h|_{C_q}$ is injective, and (iii) holds.
As in the proof of Lemma \ref{Glemma}, the statement about the generating set $G$ follows quickly.

Finally, suppose $h$ is an idempotent.  Since $\hb\in E(\T_m)$ and $\rank(\hb)=m-1$, it follows that
$\hb=e_{ab}$ for some $a,b\in[m]$ with $a\not=b$.  Since $(\gb_1\cdots\gb_{l-1})|_L=\id_L$ but $(\gb_1\cdots\gb_{l-1}\hb)|_L\not=\id_L$, we again conclude that $a\in[m]\sm L$ and $b\in L$.  We observe that $b=i$.  Indeed, if this was not the case, then we would have
$|C_bh|\leq|C_a|=n_a<n_b$, 
contradicting the fact that $h_q$ is injective for all $q\in L\sm\{i\}$.  
In particular, $\hb=e_{ai}$.  So $h$ maps $C_i$ into $C_a$, which gives $n_j=|C_ih|\leq|C_a|=n_a$.  But also $n_a\leq n_j$ since $a\in[m]\sm L$, so it follows that $n_a=n_j$.  So far we know that
\[
h=[h_1,\ldots,h_i,\ldots,h_m;e_{ai}].
\]
Since we know that $h$ maps $C_i$ into $C_a$ and $|C_ih|=n_j=|C_a|$, it follows that $h_i\in\Surj(n_i,n_j)$.  We wish to show that $h=e_{ai;h_i}$.  This will complete the proof of the lemma (with $k=a$ and $f'=h_i$).  It remains to show that $h_q=1$ for all $q\in [m]\sm\{i\}$.  In fact, since $h$ is an idempotent, and since $\im(\hb)=\im(e_{ai})=[m]\sm\{i\}$, we already know that $h_q\in E(\T_{n_q})$ for all $q\in[m]\sm\{i\}$, so it suffices to show that each such $h_q$ is injective.  We already know this is the case for $q\in L\sm\{i\}$.  We also know that $C_a=C_ih\sub C_ah$ by Proposition \ref{idempotent_prop}(iii), so it follows that $h_a$ is surjective.  But all surjective transformations of a finite set are injective, so it follows that $h_a$ is injective.
It remains to establish the injectivity of each $h_q$ with $q\in[m]\sm(L\cup\{a\})$.  To do this, we must consider two separate cases.  To simplify notation, put $u=g_1\cdots g_{l-1}$ and  $v=g_l\cdots g_r$, and write $u=[u_1,\ldots,u_m;\ub]$ and $v=[v_1,\ldots,v_m;\vb]$.  Since $e_{ji;f}=uv = [u_1v_{1\ub},\ldots,u_mv_{m\ub};\ub\;\! \vb]$, it follows that $u_q$ is injective for all $q\in[m]\sm\{i\}$.  Similarly, $u_qh_{q\ub}$ is injective for all such $q$.

{\bf Case 1.}  Suppose first that $\ub=1$.  Then $u_q\in \E_{n_q}$ for each $q\in[m]$.  Since $u_q$ is injective for each $q\in[m]\sm\{i\}$, it follows that $u_q=1$ for each such $q$.  So
\[
uh = [u_1h_1,\ldots,u_mh_m;e_{ai}] = [h_1,\ldots,h_{i-1},u_ih_i,h_{i+1},\ldots,h_m;e_{ai}].
\]
In particular, for any $q\in[m]\sm\{i\}$, $h_q=u_qh_q=u_qh_{q\ub}$ is injective, as noted above.
This completes the proof in this case.


{\bf Case 2.}  Finally, suppose $\ub\not=1$.  So $\ub\in\TmSm$.  Since $\ve_{ij}=\ker(e_{ji})=\ker(\ub\;\! \vb) \supseteq\ker(\ub)\not= \Delta$, it follows that $\ker(\ub)=\ve_{ij}$.  We claim that $a\not\in\im(\ub)$.  Indeed, suppose this was not the case, so $a=c\ub$ for some $c\in[m]$.  Since $\ub|_L=\id_L$ and $j\ub=i\ub$, it follows that $\ub$ maps $L\cup\{j\}$ into $L$, so $c\in[m]\sm(L\cup\{j\})$.  But then
\[
c e_{ji} = c\gb_1\cdots\gb_r = c \ub e_{ai}\gb_{l+1}\cdots\gb_r = ae_{ai}\gb_{l+1}\cdots\gb_r = ie_{ai}\gb_{l+1}\cdots\gb_r = i\ub e_{ai}\gb_{l+1}\cdots\gb_r = ie_{ji}=j,
\]
a contradiction, since $c\not\in\{i,j\}$.  This completes the proof of the claim that $a\not\in\im(\ub)$.  Since $\rank(\ub)=m-1$, it follows that $\im(\ub)=[m]\sm\{a\}$.
Let $Q=[m]\sm L$, and put
\[
Y=\bigcup_{q\in Q\sm\{j\}} C_q \AND Z=\bigcup_{q\in Q\sm\{a\}}C_q.
\]
Since $|C_a|=n_j=|C_j|$, it follows that $|Y|=|Z|$.  Since $\ub$ maps $Q\sm\{j\}$ bijectively onto $Q\sm\{a\}$, and since $u_p$ is injective for all $p\in Q$, it follows that $u_p$ is bijective for each $p\in Q\sm\{j\}$.  Now let $q\in Q\sm\{a\}$ be arbitrary.  Since $Q\sm\{a\}=[m]\sm(L\cup\{a\})$, the proof of the lemma will be complete if we can show that $h_q$ is injective.  Put $p=q\ub^{-1}\in Q\sm\{j\}$.  Note that
\[
uh = [u_1h_{1\ub},\ldots,u_ph_q,\ldots,u_mh_{m\ub};\ub e_{ai}].
\]
Since $u_ph_q$ is injective, and since $u_p$ is bijective, it follows that $h_q$ is injective.  As noted above, this completes the proof of the lemma. \epf


We are now able to give a lower bound for $\rank(S)$.  The next result is stated in terms of the parameters $\mu_i,\nu_i,\rho_i$ introduced at the beginning of this section and in Theorem \ref{thm:EX}.

\ms
\begin{cor}\label{cor:rank}
We have $\rank(S)\geq\rho$, where
\[
\rho = 
\sum_{i=1}^n
\left(
\mu_i\rho_i + i!{\mu_i\choose2} + \mu_i\nu_i
\right)
+ \sum_{\oijn} \mu_i\mu_j \frac{j!}{(j-i)!}.
\]
\end{cor}

\pf Let $G$ be an arbitrary generating set for $S$.  By Lemma \ref{Elemma}, $G$ contains a generating set for $\E_{n_r}^{(r)}$ for each $r\in[m]$.  These are pairwise disjoint, and each has size at least $\rank(\E_{n_r})=\rho_{n_r}$, so $G$ contains at least
\begin{equation}\tag{\ref{cor:rank}.1}\label{1.1}
\sum_{r=1}^m \rho_{n_r} = \sum_{i=1}^n \mu_i\rho_i
\end{equation}
transformations coming from these generating sets of $\E_{n_1}^{(1)},\ldots,\E_{n_m}^{(m)}$.

Next, fix some $i\in[n]$ with $\mu_i\geq2$.  Lemma \ref{Glemma} tells us that for each $p,q\in M_i$ with $p<q$, and for each $f\in\S_i$, $G$ contains some transformation $g$ such that 
\begin{itemize}\begin{multicols}3
\item[(i)] $\ker(\gb)=\ve_{pq}$,
\item[(ii)] $g_1,\ldots,g_m$ are injective, 
\item[(iii)] $g_p\in\S_i$ and $g_qg_p^{-1}=f$.
\emc\eit
(For future reference, we note that if this $g$ is an idempotent, then Lemma~\ref{Glemma} gives $g=e_{ij;f}$ or $e_{ji;f^{-1}}$.)  There are ${\mu_i\choose2}$ such $p,q$, and there are $i!$ such $f$.  Summing over all appropriate $i$, and noting that ${\mu_i\choose2}=0$ if $\mu_i\leq1$, we see that~$G$ contains at least
\begin{equation}\tag{\ref{cor:rank}.2}\label{1.2}
\sum_{i\in[n] \atop \mu_i\geq2} i!{\mu_i\choose2} = \sum_{i=1}^n i!{\mu_i\choose2}
\end{equation}
transformations of this type.  

Next, suppose $1\leq p<q\leq m$ are such that $n_p>n_q$.  Let $f\in\Inj(n_q,n_p)$ be arbitrary.  Lemma~\ref{Glemma} tells us that $G$ must contain a transformation $g$ such that 
\begin{itemize}\begin{multicols}3
\item[(i)] $\ker(\gb)=\ve_{pq}$,
\item[(ii)] $g_1,\ldots,g_m$ are injective, 
\item[(iii)] $g_p\in\S_{n_p}$ and $g_qg_p^{-1}=f$.
\emc\eit
There are $|\Inj(n_q,n_p)|=n_p!/(n_p-n_q)!$ such transformations.  For $\oijn$, there are $\mu_i\mu_j$ choices of $1\leq p<q\leq m$ with $j=n_p$ and $i=n_q$, so $G$ contains at least
\begin{equation}\tag{\ref{cor:rank}.3}\label{1.3}
\sum_{1\leq p<q\leq m \atop n_p>n_q} \frac{n_p!}{(n_p-n_q)!} = \sum_{\oijn} \mu_i\mu_j \frac{j!}{(j-i)!}
\end{equation}
transformations of this type.  

Finally, let $i\in[n]$ be such that $\mu_i\not=0$, and suppose $p\in[m]$ is such that $n_p=i$.  Let $L=\set{q\in[m]}{n_q\geq i}$.  If $1\leq j<i$ is such that $\mu_j\not=0$, then Lemma \ref{Glemma2} says that $G$ must contain some transformation $g$ such that 

~\ \ (i) $\rank(\gb)=m-1$, \qquad\qquad\qquad (ii) $|C_pg|=j$, and \qquad\qquad\qquad (iii) $g_q$ is injective for all $q\in L\sm\{i\}$.

There are $\mu_i$ such $p$, and $\nu_i$ such $j$.  So $G$ contains at least
\begin{equation}\tag{\ref{cor:rank}.4}\label{1.4}
\sum_{i\in[n] \atop \mu_i\not=0} \mu_i\nu_i = \sum_{i=1}^n \mu_i\nu_i 
\end{equation}
transformations of this type.  Finally, adding equations (\ref{1.1}--\ref{1.4}) shows that $|G|\geq\rho$.
Since $G$ is an arbitrary generating set, the result follows. \epf


We show below that this lower bound for $\rank(S)$ is precise (see Theorem \ref{thm:rank}).  In fact, we will also show that $\idrank(S)=\rank(S)$ apart from the special case in which $\mu_1=2$.  In order to deal with that case, we need Lemma~\ref{lem:U_r} below, which will also be useful when we later classify and enumerate the minimal idempotent generating sets for $S$ (Theorem \ref{thm:igs_class}).  But first we need a number of technical results.

Let $i\in[n]$.  Recall that $M_i=\set{q\in[m]}{n_q=i}$.  Let $X_i=\bigcup_{q\in M_i}C_q$, and put
\[
S_i = \set{f\in S}{f|_{X\sm X_i}=\id_{X\sm X_i},\;\! X_if\sub X_i}.
\]
The reader should not confuse $S_i$ with $\S_i$, the symmetric group on $[i]$.  Let $\P_i=\set{C_q}{q\in M_i}$.  So $\P_i$ is a uniform partition of $X_i$ into $\mu_i$ blocks of size $i$.  We aim to show that $S_i$ is isomorphic to $\EXPi$, the idempotent generated subsemigroup of $\T(X_i,\P_i)$.  The following was proved in~\mbox{\cite[Proposition 4.1]{DE1}}.

\ms
\vbox{
\begin{prop}\label{prop:EXPi}
Let $i\in[n]$ and write $M_i=\{q_1,\ldots,q_{\mu_i}\}$.  Then $f=[f_{q_1},\ldots,f_{q_{\mu_i}};\fb]\in\TXPi$ belongs to $\EXPi$ if and only if one of the following holds:
\begin{itemize}\begin{multicols}2
\itemit{i} $\fb=1$ and $f_{q_1},\ldots,f_{q_{\mu_i}}\in\E_i$,
\itemit{ii} $\fb\in\T_{X_i}\sm\S_{X_i}$. \epfres
\emc\eit
\end{prop}
}

\ms
\begin{lemma}\label{lem:Si}
Let $i\in[n]$.  Then $S_i$ is isomorphic to $\EXPi$
.
\end{lemma}

\pf There is an obvious embedding $\phi:\T(X_i,\P_i)\to\TXP$ defined, for $f\in\TXPi$, by
\[
x (f\phi) = \begin{cases}
xf &\text{if $x\in X_i$}\\
x &\text{if $x\in X\sm X_i$.}
\end{cases}
\]
So $\EXPi$ is isomorphic to its image, $T=\EXPi\phi$.  It remains to show that $S_i=T$.  Clearly, $T\sub S_i$.  Conversely, suppose $f\in S_i$, and put $g=f|_{X_i}\in\TXPi$.  Obviously, $f=g\phi$, so it suffices to prove that $g\in\EXPi$.  But this follows quickly from Lemma \ref{lem:fb=1} and Proposition~\ref{prop:EXPi}.~\epf

Next, we aim to show that any idempotent generating set for $S$ must contain a generating set for each $S_i$.  To do this, we require the next technical result.

\ms
\begin{lemma}\label{lem:fg_q}
Let $r\in[n]$ with $\mu_r\not=0$, and let $h\in S$ be such that $\rank(\hb)=m-1$, $M_r\hb\sub M_r$, $|M_r\hb|=\mu_r-1$, and
$h_q$ is injective for all $q\in[m]$.  Suppose also that $g\in E\sm S_r$ is such that $\rank(\hb\gb)=m-1$ and $(hg)_q$ is injective for all $q\in[m]$.  Then $g|_Y=\id_Y$ where $Y=X_rh$.  
In particular, $(hg)|_{X_r}=h|_{X_r}$.
\end{lemma}

\pf It is clear that $(hg)|_{X_r}=h|_{X_r}$ follows from $g|_Y=\id_Y$, so we just prove the latter.  By assumption, we have $M_r\hb=M_r\sm\{a\}$ for some $a\in M_r$.  Note that $m-1=\rank(\hb\gb)\leq\rank(\gb)\leq m$.  We now break up the proof into cases, according to whether $\rank(\gb)=m$ or $m-1$.  

{\bf Case 1.}  Suppose first that $\rank(\gb)=m$.  So $\gb=1$ and $g_q\in E(\T_{n_q})$ for all $q\in[m]$.  Let $q\in M_r\sm\{a\}$, and suppose $q=p\hb$ where $p\in M_r$.  Then $(hg)_p=h_pg_{p\hb}=h_pg_q$.  But $h_p$ and $(hg)_p$ are injective, by assumption.  Since $p\hb=q\in M_r$, it follows that $h_p\in\S_r$, so in fact $g_q=h_p^{-1}(hg)_p$ is injective.  Since also $g_q\in E(\T_r)$, as noted above, we conclude that $g_q=1$.  Since this is true for all $q\in M_r\sm\{a\}$, and since $\gb=1$, it follows that $g|_Y=\id_Y$, as desired.

{\bf Case 2.}  Suppose now that $\rank(\gb)=m-1$.  Since $g\in E$, we must have $\gb=e_{bc}$ for some $b,c\in[m]$ with $b\not=c$.  Since $\rank(\hb\gb)=\rank(\hb)=m-1$, we cannot have both $b,c\in\im(\hb)$.  We also have $g_q\in E(\T_{n_q})$ for all $q\in\im(\gb)=[m]\sm\{c\}$.  We now consider two subcases, according to whether $a$ belongs to $\im(\hb)$ or not.

{\bf Subcase 2.1.}  Suppose first that $a\not\in\im(\hb)$.  Since $b,c$ do not both belong to $\im(\hb)$, we must have either $b=a$ or $c=a$.  Suppose first that $c=a$.  Note that $\gb|_{M_r\sm\{a\}}=\id_{M_r\sm\{a\}}$ and $g_q\in E(\T_r)$ for all $q\in M_r\sm\{a\}$ so, as in Case~1, we conclude that $g_q=1$ for all such $q$, and therefore, $g|_Y=\id_Y$.  Now suppose $b=a$.  If $c\not\in M_r$, then again $\gb|_{M_r\sm\{a\}}=\id_{M_r\sm\{a\}}$ and $g_q\in E(\T_r)$ for all $q\in M_r\sm\{a\}$, and the proof concludes as above.  So suppose $c\in M_r$.  In fact, we will show that this case is not possible.  
Since $\im(\hb)=[m]\sm\{a\}$ and $M_r\hb\sub M_r$, it follows that $\hb$ maps $[m]\sm M_r$ bijectively into $[m]\sm M_r$.  But $h_q$ is injective for all $q\in[m]\sm M_r$, so it follows that $h|_{X\sm X_r}\in\T_{X\sm X_r}$ is injective, and hence bijective.  We deduce that $h_q$ is bijective for all $q\in[m]\sm M_r$.  Let $q\in[m]\sm M_r$ and put $p=q\hb^{-1}$.  Since $(gh)_p=h_pg_q$ is injective, it follows that $g_q$ is injective.
%
%
But also $g_q\in E(\T_{n_q})$ for all such $q$, giving $g_q=1$ for $q\in[m]\sm M_r$.  Since also $\gb|_{[m]\sm M_r}=\id_{[m]\sm M_r}$, we would have $g\in S_r$, a contradiction.  This completes the proof in Subcase 2.1.

{\bf Subcase 2.2.}  Finally, suppose $a\in\im(\hb)$.  
As before, if $c\not\in M_r\sm\{a\}$, then $g|_Y=\id_Y$ quickly follows.  So suppose $c\in M_r\sm\{a\}$.  Actually, we will show that this is impossible.  In particular, $c\in\im(\hb)$, so $b\not\in\im(\hb)$.  Next we claim that $n_b<r$.  Indeed, suppose this is not the case.  Note first that $b\not\in\im(\hb)$ implies $n_b\not=r$, since $M_r\sub\im(\hb)$.  So we must have $n_b>r$.  Then some $q\in[m]$ with $n_q>r$ is mapped by $\hb$ to some $p\in[m]$ with $n_p\leq r<n_q$.  But then $|C_qh|\leq|C_p|<|C_q|$, contradicting the fact that $h_q$ is injective.  This completes the proof of the claim that $n_b<r$.  But now let $q\in M_r$ be such that $q\hb=c$.  Then $|C_qhg|\leq|C_cg|\leq|C_b|=n_b<r=|C_q|$, contradicting the assumption that $(hg)_q$ is injective.  This completes the proof. \epf


\begin{lemma}\label{lem:U_r}
Let $U$ be an arbitrary idempotent generating set for $S$ and let $r\in[n]$.  Then $U\cap S_r$ is a generating set for $S_r$.
\end{lemma}

\pf If $\mu_r=0$, then $S_r=\{1\}=\la\emptyset\ra$, and the result is trivially true.  So suppose $\mu_r\geq1$, and write $U_r=U\cap S_r$.  By Lemma \ref{lem:Si} and \cite[Corollary 4.2]{DE1}, $S_r$ is generated by all elements of the form
\ms\ms\ms\bit\bmc2
\item[(a)] $e_{ij}^{(k)}$ for $i,j\in[r]$ with $i\not=j$ and all $k\in M_r$, and
\item[(b)] $e_{ij;f}$ for $i,j\in M_r$ with $i\not=j$ and all $f\in\S_r$.
\emc\eit
So it suffices to show that $\la U_r\ra$ contains each element of type (a) and (b).
Now, for each $k\in M_r$,~$U$ contains a generating set $V$ for $\E_{n_k}^{(k)}=\E_r^{(k)}$ by Lemma \ref{Elemma}, and we clearly have $V\sub U_r$.  In particular, $\la U_r\ra$ contains all elements of type (a).  So now fix $f\in\S_r$ and $i,j\in M_r$ with $i\not=j$.  Consider an expression
$
e_{ij;f} = g_1\cdots g_k,
$
where $g_1,\ldots,g_k\in U\sm\{1\}$.  Let $L=\set{q\in[k]}{g_q\in U_r}$, and write $L=\{q_1,\ldots,q_l\}$ where $q_1<\cdots<q_l$.  We first aim to show that $(g_1\cdots g_k)|_{X_r}=(g_{q_1}\cdots g_{q_l})|_{X_r}$.  

For $p\in[k]$, let $h_p=g_1\cdots g_p$, and write $h_p=[h_{p1},\ldots,h_{pm};\hb_p]$.  We claim that for all $p\in[k]$,
\begin{itemize}
\begin{multicols}2
\item[(i)] $h_{pq}$ is injective for all $q\in[m]$,
\item[(ii)] $\rank(\hb_p)=m-1$,
\item[(iii)] $M_{r}\hb_p\sub M_{r}$,
\item[(iv)] $|M_{r}\hb_p|=\mu_{r}-1$.
\emc\eit
By Lemma \ref{Glemma}, $h_1=g_1=e_{ij;f}$ or $e_{ji;f^{-1}}$, so the claim is true for $p=1$.  Now suppose (i--iv) all hold for some $1\leq p<k$.  Since $e_{ij;f}=h_{p+1}g_{p+2}\cdots g_k$, it is clear that each $h_{p+1,q}$ is injective, so (i) holds.  Next, note that $m-1=\rank(\gb_1)\leq\rank(\hb_{p+1})\leq\rank(e_{ij})=m-1$, giving (ii).  Next, we have $|M_{r}\hb_{p+1}|\leq|M_{r}\gb_1|=\mu_{r}-1$.  But $\rank(\hb_{p+1})=m-1$, so $|A\hb_{p+1}|\geq|A|-1$ for any subset $A\sub[m]$, and (iv) follows.  For (iii), note that the induction hypothesis gives $M_{r}\hb_p\sub M_{r}$.  If $g_{p+1}\in S_{r}$, then $M_{r}\gb_{p+1}\sub M_{r}$ by definition, so $M_{r}\hb_{p+1}=M_{r}\hb_p\gb_{p+1}\sub M_{r}$.  If $g_{p+1}\not\in S_{r}$, then all the conditions of Lemma \ref{lem:fg_q} are satisfied (with $h=h_p$ and $g=g_{p+1}$).  We conclude then that $g_{p+1}|_Y=\id_Y$, where $Y=X_{r}h_p$.  In particular, $\gb_{p+1}|_{M_r\hb_p}=\id_{M_r\hb_p}$.  But then $M_{r}\hb_{p+1}=M_{r}\hb_p\gb_{p+1} = M_{r}\hb_p\sub M_{r}$.  This completes the proof of the inductive step and, hence, of the claim.

Now suppose $p\in[k]$ is such that $g_p\not\in S_r$.  In particular, $p\geq2$ (since $g_1\in S_r$, as noted above).  By the above claim (and as noted in its proof), the conditions of Lemma \ref{lem:fg_q} are satisfied (for $h=h_{p-1}$ and $g=g_p$), so we conclude that $h_p|_{X_r}=h_{p-1}|_{X_r}$.  So, if $Q=[k]\sm L=\set{p\in[k]}{g_p\not\in U_r}$, and $Q=\{p_1,\ldots,p_{s}\}$ where $s=k-l$ and $p_1<\cdots<p_{s}$, then
\begin{align*}
e_{ij;f}|_{X_r} = (g_1\cdots g_k)|_{X_r} &= h_{p_{s}}|_{X_r} (g_{p_{s}+1}\cdots g_k) \\
&= h_{p_{s}-1}|_{X_r} (g_{p_{s}+1}\cdots g_k) \\
&= h_{p_{s-1}}|_{X_r} (g_{p_{s-1}+1}\cdots g_{p_{s}-1}) (g_{p_{s}+1}\cdots g_k) \\
&= h_{p_{s-1}-1}|_{X_r} (g_{p_{s-1}+1}\cdots g_{p_{s}-1}) (g_{p_{s}+1}\cdots g_k) \\
&\hspace{1.9mm} \vdots \\
&= (g_1\cdots g_{p_1-1})|_{X_r} (g_{p_1+1}\cdots g_{p_2-1}) \cdots (g_{p_{s}+1}\cdots g_k)\\
&= ((g_1\cdots g_{p_1-1}) (g_{p_1+1}\cdots g_{p_2-1}) \cdots (g_{p_{s}+1}\cdots g_k))|_{X_r} \\
&= (g_{q_1}\cdots g_{q_l})|_{X_r}.
\end{align*}
But also $e_{ij;f}|_{X\sm X_r} = \id_{X\sm X_r} = (g_{q_1}\cdots g_{q_l})|_{X\sm X_r}$, so it follows that $e_{ij;f}=g_{q_1}\cdots g_{q_l}\in\la U_r\ra$, completing the proof. \epf

\begin{cor}\label{cor:idrank}
If $\mu_1=2$, then $\idrank(S)\geq\rho+1$, where $\rho$ is defined in Corollary \ref{cor:rank}.
\end{cor}

\pf Let $G$ be an arbitrary idempotent generating set for $S$.  By Lemma \ref{lem:U_r}, $G$ contains a generating set $U_i$ for $S_i\cong\EXPi$ for each $i\in[n]$, and we have $|U_i|\geq\rank(\EXPi)=\rho_{\mu_i,i}$.  By Theorem \ref{thm:rank_uniform}, $\rho_{\mu_i,i}=\mu_i\rho_i + i!{\mu_i\choose2}$, unless $i=1$, in which case $\rho_{\mu_i,i}=\rho_{21}=2=1+\left(\mu_1\rho_1 + 1!{\mu_1\choose2}\right)$.
So $G$ contains at least
\[
\sum_{i=1}^n \rho_{\mu_i,i} 
= 1 + \sum_{i=1}^n \left( \mu_i\rho_i + i!{\mu_i\choose2} \right)
\]
elements from these generating sets of $S_1,\ldots,S_n$.  By the last two paragraphs of the proof of Corollary~\ref{cor:rank},~$G$ contains an additional
\[ 
\sum_{i=1}^n
 \mu_i\nu_i
+ \sum_{\oijn} \mu_i\mu_j \frac{j!}{(j-i)!}
\]
elements.  Adding the two expressions above, we see that $|G|\geq\rho+1$, and the proof is complete. \epf

For the proof of the next result, we use the standard notation $f=\trans{A_1 & \cdots & A_k}{a_1 & \cdots & a_k}$ to indicate that $f$ is the function with domain $A_1\cup\cdots\cup A_k$ that maps each of the points in $A_q$ to $a_q$ for each $q\in[k]$.

\ms
\begin{prop}\label{prop:U}
For each $q\in[n]$, let $U_q$ be an idempotent generating set for $S_q$. 
For each $i\in[m]$, choose sets $J_i\sub[m]$ such that $|J_i|=\nu_{n_i}$ and $\set{n_j}{j\in J_i}=\set{n_q}{q\in[m],\;\! n_q<n_i}$.  For each $i\in[m]$ and $j\in J_i$, choose some $f_{ij}\in\Surj(n_i,n_j)$.
Then $U=U_1\cup\cdots\cup U_n \cup W_1\cup W_2$ is a generating set for $S$, where
\[
W_1 = \set{e_{ij;f}}{\oijm,\;\! n_i>n_j,\;\! f\in\Inj(n_j,n_i)}  \AND W_2 = \set{e_{ji;f_{ij}}}{i\in[m],\;\! j\in J_i}.
\]
\end{prop}

\pf By Proposition \ref{prop:G1G2}, it suffices to prove that $G_1\cup G_2\sub \la U\ra$.  Since $S_q=\la U_q\ra$ for each $q\in[n]$, it follows that~$\la U\ra$ contains
\bit
\item[(i)] each $e_{ij}^{(k)}$ with $k\in[m]$, $i,j\in[n_k]$ and $i\not=j$, and
\item[(ii)] each $e_{ij;f}$ with $i,j\in[m]$, $i\not=j$, $n_i=n_j$ and $f\in\S_{n_i}$.
\eit
Since $W_1\sub U$, it remains to show that $\la U\ra$ contains
\bit
\item[(iii)] each $e_{ji;f}$ with $i,j\in[m]$, $n_i>n_j$ and $f\in\Surj(n_i,n_j)$.
\eit
Let $i,j,f$ be as in (iii).  
Let $k\in J_i$ be such that $n_k=n_j$, and for simplicity, put $g=f_{ik}$.  So $e_{ki;g}\in U$.  Write $f=\trans{A_1 & \cdots & A_{n_j}}{1 & \cdots & n_j}$ and $g=\trans{B_1 & \cdots & B_{n_j}}{1 & \cdots & n_j}$.  Also, choose $b_1,\ldots,b_{n_j}$ such that $b_q\in B_q$ for each $q$.  Put $h=\trans{A_1 & \cdots & A_{n_j}}{b_1 & \cdots & b_{n_j}}\in\T_{n_i}$.  Since $n_j<n_i$, note that $h\in\E_{n_i}$.  So $h^{(i)}\in\E_{n_i}^{(i)}\sub\la U\ra$.  It is clear that
$
e_{ki;f} = h^{(i)} e_{ki;g} \in\la U\ra.
$
In particular, if $k=j$, then we have shown that $e_{ji;f}=e_{ki;f}\in\la U\ra$.  So suppose $k\not=j$.  Now choose $a_1,\ldots,a_{n_j}$ so that $a_q\in A_q$ for each $q$.  Put $d=\trans{1&\cdots&n_j}{a_1&\cdots&a_{n_j}}\in\Inj(n_j,n_i)$.  Note that $e_{ij;d}\in W_1\sub U$, and that $e_{jk;1}\in\la U\ra$ as shown above, since $n_j=n_k$.  It is clear that $df=1\in\S_{n_j}$, and one may then easily check that
$
e_{ji;f} = (e_{ij;d} e_{jk;1} e_{ki;f})^2\in\la U\ra,
$ 
completing the proof. \epf

We are now ready to prove the two main results of the paper.  Again, these are stated in terms of the parameters $\mu_i,\nu_i,\rho_i$ introduced at the beginning of this section and in Theorem \ref{thm:EX}. 

\ms
\begin{thm}\label{thm:rank}
We have $\rank(S)=\idrank(S)=\rho$, where
\[
\rho=
\sum_{i=1}^n
\left(
\mu_i\rho_i + i!{\mu_i\choose2} + \mu_i\nu_i
\right)
+ \sum_{\oijn} \mu_i\mu_j \frac{j!}{(j-i)!},
\]
unless $\mu_1=2$ in which case $\rank(S)=\rho$ and $\idrank(S)=\rho+1$.
\end{thm}

\pf Let $U=U_1\cup\cdots\cup U_n \cup W_1\cup W_2$ be an idempotent generating set as described in Proposition \ref{prop:U}, with $U_1,\ldots,U_n$ of minimal size.  As in the proof of Corollary \ref{cor:rank}, 
\[
|W_1| = \sum_{\oijn} \mu_i\mu_j \frac{j!}{(j-i)!}
\AND
|W_2| = \sum_{r=1}^m\nu_{n_r} = \sum_{i=1}^n\mu_i\nu_i.
\]
As in the proof of Corollary \ref{cor:idrank}, 
\[
|U_1|+\cdots+|U_n| = \sum_{i=1}^n \left( \mu_i\rho_i + i!{\mu_i\choose2} \right),
\]
unless $\mu_1=2$, in which case we must add $1$ to the right hand side of this last expression.  Adding these values, we conclude that
\[
|U| = \begin{cases}
\rho &\text{if $\mu_1\not=2$}\\
\rho+1 &\text{if $\mu_1=2$.}
\end{cases}
\]
Combined with Corollaries \ref{cor:rank} and \ref{cor:idrank}, and noting that $U\sub E$, this shows that $\rank(S)=\idrank(S)=\rho$ if $\mu_1\not=2$, and also that $\idrank(S)=\rho+1$ if $\mu_1=2$.  To complete the proof, it suffices to prove that $S=\la V\ra$ for some $V\sub S$ with $|S|=\rho$ if $\mu_1=2$.  For the remainder of the proof, we assume that $\mu_1=2$.

We have already seen that $S=\la U\ra$ and $|U|=\rho+1$.  Also, since there is a unique generating set of size $2$ for $S_1\cong\E(X_1,\P_1)\cong\E_2$, namely $U_1=\{f,g\}$ where $f=e_{m-1,m;1}$ and $g=e_{m,m-1;1}$, we see that $U$ must contain both $f$ and $g$.  Let $h\in\Inj(1,n_1)=\Inj(n_m,n_1)$ be arbitrary, and put $e=e_{1m;h}$.  So $e\in W_1\sub U$.  It is easy to check that $e=(eg)f$ and $g=f(eg)$.  It follows that $\la e,f,g\ra=\la eg,f\ra$ and so $S=\la V\ra$, where $V=(U\sm\{e,g\})\cup\{eg\}$.  Since $|V|=\rho$, this completes the proof.  \epf 

For the statement of the next result, by ``minimal idempotent generating set'' we mean an idempotent generating set that has the smallest possible size.

\ms\ms
\begin{thm}\label{thm:igs_class}
\bit
\itemit{i} Every idempotent generating set of $S$ contains a minimal idempotent generating set.
\itemit{ii} Every minimal idempotent generating set of $S$ is of the form described in Proposition \ref{prop:U} and with each $U_q$ of minimal size.
\itemit{iii} The number of minimal idempotent generating sets of $S$ is equal to
\[
\prod_{i=1}^n \si_{\mu_i,i} \times \prod_{\oijn \atop \mu_i\not=0\not=\mu_j} \mu_i\mu_j S(j,i)i!,
\]
where $S(j,i)$ is a Stirling number (of the second kind) and the numbers $\si_{\mu_i,i}$ are defined in Theorem~\ref{thm:enum_uniform}.  
\eit
\end{thm}

\pf Let $U$ be an arbitrary idempotent generating set for $S$.  By Lemma \ref{lem:U_r}, $U$ contains an idempotent generating set $U_r$ of $S_r$ for each $r\in[n]$.  By Theorem \ref{thm:enum_uniform}, each $U_r$ contains an idempotent generating set $V_r$ of minimal size.  As in the proof of Corollary \ref{cor:rank}, $U$ must contain the sets
\begin{align*}
W_1 = \set{e_{ij;f}}{\oijm,\;\! n_i>n_j,\;\! f\in\Inj(n_j,n_i)} \AND
W_2 = \set{e_{ji;f_{ij}}}{i\in[m],\;\! j\in J_i},
\end{align*}
for some choice of sets $J_i$ and  functions $f_{ij}\in\Surj(n_i,n_j)$.  The set $V_1\cup\cdots\cup V_n\cup W_1\cup W_2\sub U$ has size $\idrank(S)$, as stated in Theorem \ref{thm:rank}, and is a generating set for $S$ by Proposition \ref{prop:U}.  This completes the proof of (i).

If $U$ is an arbitrary idempotent generating set of minimal possible size, then we must in fact have $U=V_1\cup\cdots\cup V_n\cup W_1\cup W_2$ (in the above notation), proving (ii).  For each $i\in[n]$, we may choose $V_i$ in $\si_{\mu_i,i}$ ways.  To specify $W_2$, for each $k\in M_j$ where $j\in[n]$ is such that $\mu_j\not=0$, and for each $i\in[n]$ with $i<j$ and $\mu_i\not=0$, we must choose some $e_{kl;f}$ where $l\in M_i$ and $f\in\Surj(j,i)$.  There are $\mu_j$ such $k$, $\mu_i$ such $l$, and $|\Surj(j,i)|=S(j,i)i!$ such~$f$.  Multiplying these values as appropriate gives (iii).~\epf

\footnotesize \def\bibspacing{-1.1pt} \bibliography{endp2_biblio} \bibliographystyle{plain} \end{document}